\documentclass[11pt]{amsart}
\usepackage{amssymb,amsmath,amsthm,latexsym,cancel}
\usepackage{graphicx}
\usepackage{epsfig}
\newtheorem{theorem}{Theorem}[section]
\newtheorem{lemma}[theorem]{Lemma}

\newtheorem{corollary}[theorem]{Corollary}
\newtheorem{problem}[theorem]{Problem}
\theoremstyle{remark}
\newtheorem{definition}[theorem]{Definition}
\newtheorem{remark}[theorem]{Remark}
\newtheorem{example}[theorem]{Example}

\newcommand\R{{\mathcal R}}

\newcommand\Z{{\mathbb Z}}
\newcommand\C{\mathfrak{C}}
\newcommand\D{\mathfrak{D}}
\newcommand\E{\mathfrak{E}}
\newcommand\K{\mathfrak{K}}

\begin{document}

\title{Categories of diagrams with irreversible moves}
\date{February 6, 2013}
\author{Maciej Niebrzydowski}
\address[Maciej Niebrzydowski]{Institute of Mathematics\\ 
Faculty of Mathematics, Physics and Informatics\\
University of Gda{\'n}sk, 80-308 Gda{\'n}sk, Poland}
\address{Department of Mathematics\\
	 University of Louisiana at Lafayette\\
	 Lafayette, LA 70504-1010, USA}
\email{mniebrz@gmail.com}


\thispagestyle{empty}

\begin{abstract}
We work with a generalization of knot theory, in which one diagram is reachable from another via a finite sequence of moves if a fixed condition regarding the existence of certain morphisms in an associated category is satisfied for every move of the sequence. This conditional setting leads to a possibility of irreversible moves, terminal states, and to using functors more general than the ones used as knot invariants. Our main focus is the category of diagrams with a binary relation on the set of arcs, indicating which arc can move over another arc. We define homology of binary relations, and merge it with quandle homology, to obtain the homology for partial quandles with a binary relation. This last homology can be used to analyze link diagrams with a binary relation on the set of components.  
\end{abstract}

\maketitle

\section{Diagrams and categories}
In this paper, we propose to study conditional diagram theories, where an elementary move on a diagram is allowed if a fixed condition, expressed in terms of existence of certain morphisms, is satisfied. We are mostly interested in transforming a diagram theory that has a Reidemeister-type theorem 
\footnote{That is, two diagrams $D_1$ and $D_2$ represent the same object iff there is a finite sequence of elementary reversible moves leading from $D_1$ to $D_2$.} into conditional diagram theories.

For a given category $\C$, denote its objects by $\C_0$, and its morphisms by $\C_1$. 

Let $\D$ be a small category \footnote{A category is small if both the collection of its objects and the collection of its morphisms are sets.} whose objects $\D_0$ are diagrams, possibly taken up to some equivalence, and morphisms $\D_1$ are sequences of reversible elementary moves on the diagrams. We will call such $\D$ a category of diagrams and moves. Examples include: diagrams of classical knots with sequences of Reidemeister moves, diagrams of virtual knots with sequences of virtual Reidemeister moves, diagrams of spatial graphs with sequences of graphical Reidemeister moves, and diagrams of knotted surfaces with sequences of Roseman moves; in the first three cases diagrams are placed on the plane and considered up to planar isotopy, in the last case they are taken up to ambient isotopy in the $3$-space in which they are located. 

An invariant in such a category is a function $h$ defined on diagrams and assigning to them objects of some category $\C$, in such a way that if a diagram $D_2$ is obtained from a diagram $D_1$ by an elementary move of a type permitted by the theory (e.g., Reidemeister move) \footnote{We assume that the existence of an elementary move between the diagrams $D_1$ and $D_2$ implies that they are not the same as objects of the category $\D$.}, then $h(D_1)$ and $h(D_2)$ are isomorphic as objects of the category $\C$. Thus, we can think of a pair of objects $(D_1, D_2)$, and the corresponding pair $(C_1, C_2)=(h(D_1), h(D_2))$, such that there is an elementary move $\alpha\colon D_1\to D_2$ and an isomorphism $\beta\colon C_1\to C_2$. 

This situation can be viewed from a different angle: suppose that the move $\alpha$ {\it is permitted because} there is an isomorphism $\beta$ between the objects of the companion pair $(C_1, C_2)=(h(D_1), h(D_2))$, for a given function $h$. If $h$ is an invariant, this point of view does not give anything new,
but once we allow some more general ways of assigning pairs of objects of $\C$ to the pair of diagrams $(D_1, D_2)$, and consider more general morphisms in place of isomorphisms, we obtain a nontrivial generalization of diagram theories.

The conditional setting naturally leads to the use of (elementary) logic. In this introductory paper, we will assume that the number of pairs of objects assigned to $(D_1, D_2)$ is finite, and that propositional calculus is used, as this is sufficient for all the examples that we consider; more generally, a higher order logic could be used.

Recall that in propositional calculus one starts with an infinite set of atomic formulas (also called propositional variables) which take values in the set $\{$ True, False $\}$ (denoted in this paper by $\{T,F\}$ or $\{1,0\}$). Admissible compound statements are formed using logical connectives:  $\vee$ stands for `or', $\wedge$ for `and', and $\neg$ for `not'. The set of well-formed formulas (wffs) is defined by the rules:
\begin{itemize}
\item[1)] any atomic formula is a wff (we can also include T and F as wffs);
\item[2)] if $w$ and $v$ are wffs, then so are $w\vee v$, $w\wedge v$, and $\neg w$;
\item[3)] any wff is created via a finite number of applications of 1) and 2).
\end{itemize} 
By assigning a value from $\{T,F\}$ to each atomic formula, the value from $\{T,F\}$ is assigned to each wff.

The general Definition \ref{genknot} that we will now give will be followed by some concrete examples of conditional categories. The idea is to take a category $\D$ of diagrams and moves, and use the set of its objects $\D_0$, and a subset of its morphisms satisfying a fixed condition, to create a new category.

\begin{definition} \label{genknot}
For a countable number of pairs of symbols of objects, 
\[(Obj_1,Obj_1'), (Obj_2,Obj_2'),\ldots,(Obj_i,Obj_i'),\ldots\] of some (yet unspecified) category, take the set of atomic formulas $a_1$, $a_2,\ldots$, $a_i,\ldots$, defined as:\\ 
$a_1:\ (\exists m_1\colon Obj_1\to Obj_1')$,\\
$a_2:\ (\exists m_2\colon Obj_2\to Obj_2'),\ldots$,\\
$a_i:\ (\exists m_i\colon Obj_i\to Obj_i'),\ldots$,\\
where $m_i$ denotes a morphism, for each $i=1,2,\ldots$.
Consider the set $W$ of well-formed formulas of the propositional calculus on this set of variables. Now, fix a category $\C$, and a category $\D$ of diagrams and moves. Let $w\in W$ be a well-formed formula built from the atomic formulas $a_1, a_2,\ldots,a_n$.
Let $\mathbf{D}$ consist of all pairs $(D_1,D_2)$ of diagrams from $\D_0\times\D_0$ such that there exists an elementary move $\alpha\in\D_1$, $\alpha\colon D_1\to D_2$. Let $f$ be some mapping assigning to each pair $(D_1, D_2)\in\mathbf{D}$ an $n$-element sequence of pairs of objects of $\C$, that is, \[f(D_1,D_2)=\{(C_i,C'_i)\}_{i\in \{1,\ldots,n\}},\] where $C_i$, $C'_i \in \C_0$, for $i=1,\ldots,n$. 
The $i$-th element of this sequence is now used to determine the value of the atomic formula $a_i$, for $i=1,2,\ldots n$, by taking $Obj_i$ to be the symbol for $C_i$ and letting $Obj'_i$ denote $C'_i$. In this way the value $F_w(a_1,\ldots,a_n)\in\{T,F\}$ is assigned to the formula $w$. A new category, denoted by $\D\cap_f^w\C$, is formed. Its objects are the same as the objects of $\D$. Its morphisms are generated by the subset $M$ of morphisms of $\D$ determined as follows. Let $\alpha\colon D_1\to D_2$ be an elementary move of $\D$.
Then $\alpha$ is included in $M$ if and only if 
\[F_w(a_1,\ldots,a_n)=T,\ \textrm{for}\ f(D_1,D_2)=\{(C_i,C'_i)\}_{i\in \{1,\ldots,n\}}.\] 
$\D_0$ together with $M$ form a directed graph, and $\D\cap_f^w\C$ is defined to be the free category on this graph (all the compositions of morphisms and identity morphisms necessary to obtain a category are added at this stage). We call $\D\cap_f^w\C$ a {\it conditional diagram (knot, graph, etc.) theory}.
\end{definition}

\begin{remark}
If $w$ in the above definition is a tautology, then $\D\cap_f^w\C=\D$.
\end{remark}

The following is the main basic question that appears when considering conditional diagram theories.

\begin{problem}[Reachability problem]
Given diagrams $D_1$ and $D_2$ from $(\D\cap_f^w\C)_0$, is it possible to reach $D_2$ via a finite sequence of elementary moves $\alpha\in M$ allowed in $\D\cap_f^w\C$, starting from $D_1$? In other words, is there a morphism from $D_1$ to $D_2$ in $\D\cap_f^w\C$?
\end{problem}
 We write $D_1\rightarrow D_2$ if $D_2$ is reachable from 
$D_1$, and $D_1\leftrightarrow D_2$ if also $D_2\rightarrow D_1$.

\begin{definition}
An {\it indicator} is a condition that is satisfied if \mbox{$D_1\rightarrow D_2$}, for any $D_1$, $D_2\in\D\cap_f^w\C$.
\end{definition}

\begin{example} 
If $\mathcal{F}$ is a functor from $\D\cap_f^w\C$, then the existence of a morphism between
$\mathcal{F}(D_1)$ and $\mathcal{F}(D_2)$ is an indicator.
\end{example}

\begin{definition}\label{knot}
Let $\D$ be a category of diagrams and moves, where any two diagrams connected by a finite sequence of elementary reversible moves describe the same object (e.g., knot, graph, etc.). For a conditional diagram theory
$\D\cap_f^w\C$, we define the {\it object (i.e., knot, graph, etc.) of a diagram} $D\in(\D\cap_f^w\C)_0$ to be the maximal connected (both properties as an underlying unoriented subgraph) subcategory $\K(D)$ of $\D\cap_f^w\C$ such that $D$ is an object of $\K(D)$.
\end{definition} 

There are simpler notions, useful in determining $\K(D)$.
\begin{definition}
For $D$ as in Definition \ref{knot}, we define the {\it out-object of $D$ (out-knot, out-graph, etc.)}, denoted by $out(D)$, to be the full subcategory \footnote{ A subcategory $\mathcal{S}$ of a category $\mathcal{C}$ is called  full if it includes all morphisms of $\mathcal{C}$ between objects of $\mathcal{S}$.} of 
$\D\cap_f^w\C$ with the set of objects consisting of $D$, and all the diagrams $D'$ reachable from $D$. 

Similarly, let the {\it in-object of $D$}, $in(D)$, be the full subcategory of 
$\D\cap_f^w\C$ with the set of objects consisting of $D$, and all the diagrams $D'$ such that $D$ is reachable from $D'$.
\end{definition}

\begin{remark}
The relation of reachability in an oriented graph gives a preorder on the set of vertices.
It follows that a poset can be obtained by identifying diagrams $D_1$ and $D_2$ from $\D\cap_f^w\C$ such that $D_1\leftrightarrow D_2$.
\end{remark}

\begin{example}\label{firstex}
Let $\C$ be a discrete, two-element category \footnote{A discrete category is a category whose only morphisms are the identity morphisms.} with objects denoted by 0 and 1. Let $\D$
be a category of classical link diagrams with sequences of Reidemeister moves. 
Let $\hat{f}\colon \D_0\to\C$ be a function assigning to a diagram the number of its crossings modulo 2. Let $f\colon\mathbf{D}\to (\C\times\C)_0$ be defined by: \[f(D_1,D_2)=\{(\hat{f}(D_1),\hat{f}(D_2))\},\] and let $w$ be the atomic formula $(\exists m_1\colon Obj_1\to Obj_1')$. Then $\D\cap_f^w\C$ is a conditional knot theory in which the first Reidemeister move is not allowed (there is no morphism between 0 and 1, or between 1 and 0, and the first Reidemeister move is the only one that changes the parity of the number of crossings).
\end{example}

\begin{example}\label{secondex}
Let $\D$ be as in Example \ref{firstex}, and let $\hat{f}\colon\D_0\to\mathbb{R}$ be a function assigning a numerical value to every diagram. Let $\C=\mathbb{R}$, considered as a poset with the standard relation `less than or equal'. Let $f(D_1,D_2)=\{(\hat{f}(D_1),\hat{f}(D_2))\}$, and $w$ be the formula $(\exists m_1\colon Obj_1\to Obj_1')$. Then $\D\cap_f^w\C$ is a conditional knot theory in which the only Reidemeister moves allowed are the ones that do not decrease the numerical value assigned to diagrams by the function $\hat{f}$.

If $f$ is given by \[f(D_1,D_2)=\{(\hat{f}(D_1),\hat{f}(D_2)),(\hat{f}(D_2),\hat{f}(D_1))\},\] and 
\[w=(\exists m_1\colon Obj_1\to Obj_1')\wedge\neg(\exists m_2\colon Obj_2\to Obj_2'),\] then in the conditional knot theory $\D\cap_f^w\C$, a Reidemeister move between $D_1$ and $D_2$ is permitted if it increases the numerical value, i.e., $\hat{f}(D_1)<\hat{f}(D_2)$. In particular, if $\hat{f}$ assigns the number of crossings to a diagram, then the third Reidemeister move is ruled out.
\end{example}

Many examples of conditional theories can be obtained from categories $\D$ that have diagrams decorated by some categorical objects. By decorated we mean that an object $D^E\in\D_0$ consist of the underlying diagram $D$ and a set $E$ of objects of some fixed category $\E$.
As an example, consider knot or link diagrams with arcs labeled by the elements of a lattice $\mathfrak{L}$. Then, the definition of elementary Reidemeister-type moves (i.e., moves that are just the usual Reidemeister moves if we forget about the labels on arcs) in the category $\D$ can involve changes of the labels. Note that these labeled moves are in general no longer strictly `local' moves, because the labels on the arcs can extend well beyond the local pictures on which the moves are shown, and influence the set of possible moves in different regions of the diagram.

When defining the Reidemeister-type moves on labeled diagrams (that is, diagrams with labels on arcs), one always needs to consider the possibility that some of the arcs involved in the move could be the same (if connected outside of the local picture), and thus should have the same label after the move. Let us consider these requirements in more depth, using as an example the labeling of arcs by elements of a lattice $\mathfrak{L}$.

\begin{figure}
\begin{center}
\includegraphics[height=3 cm]{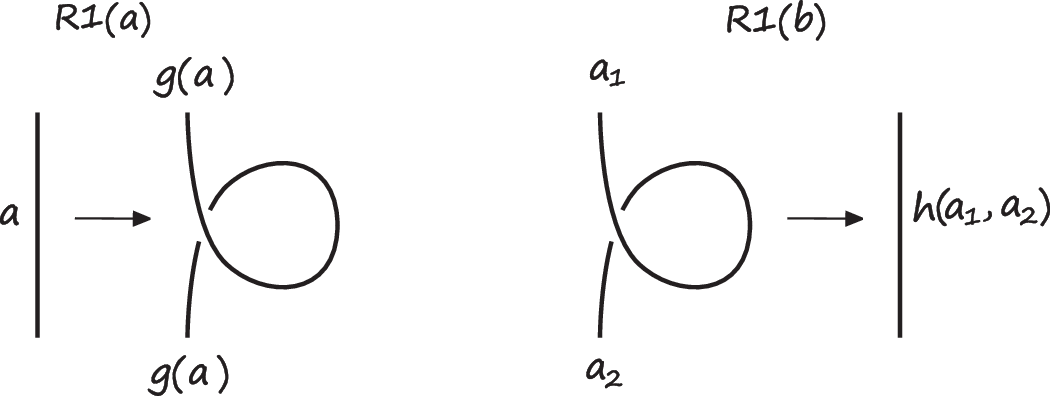}
\caption{}\label{lat1}
\end{center}
\end{figure}

In the labeled Reidemeister move of type 1(a) (see Fig. \ref{lat1}), the requirement is that after the move both arcs receive the same label; in R1(b) there are no conditions. Here, $g\colon \mathfrak{L}\to\mathfrak{L}$, and 
$h\colon \mathfrak{L}\times\mathfrak{L}\to\mathfrak{L}$ denote arbitrary functions.

\begin{figure}
\begin{center}
\includegraphics[height=3.5cm]{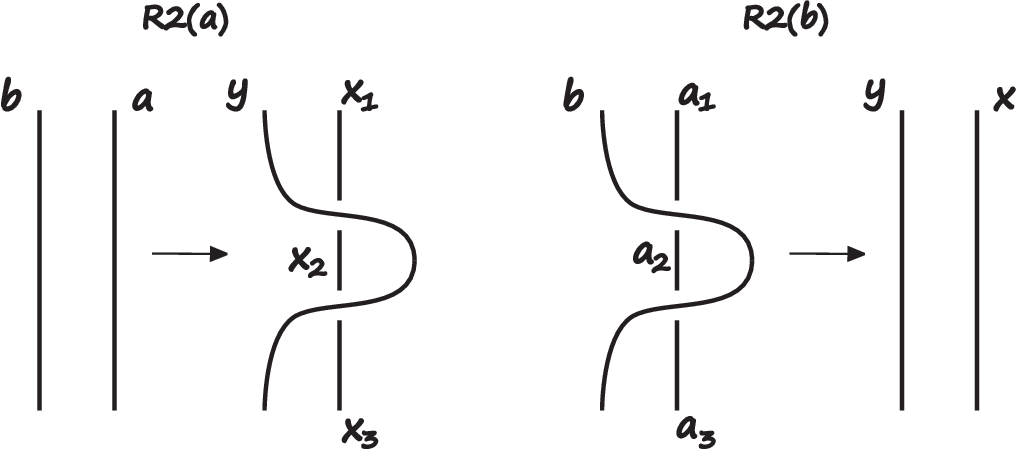}
\caption{}\label{lat2}
\end{center}
\end{figure}

In the move R2(a) illustrated in Fig. \ref{lat2}, we use 
\[f_1,\ f_2,\ f_3\colon \mathfrak{L}\times \mathfrak{L}\to \mathfrak{L}.\] Let $x_1=x_3=f_1(a,b)$, $y=f_2(a,b)$,
and $x_2=f_3(a,b)$. The function $f_3$ can be arbitrary, because the arc labeled $x_2$ is not connected to any other arc. On the other hand, it's possible that the arcs labeled $a$ and $b$ are connected, and thus it's necessary that 
\[f_1(b,b)=f_2(b,b),\ \textrm{for every}\ b\in \mathfrak{L}.\] A natural example of a pair of functions satisfying this last condition is:
\[f_1(x,y)=x\vee y\ \textrm{and}\ f_2(x,y)=x\wedge y,\]
where $\vee$ is the join, and $\wedge$ is the meet of the lattice.

For the move R2(b), consider the functions $f_1$, $f_2\colon \mathfrak{L}^4\to \mathfrak{L}$, and let 
\[x=f_1(a_1,a_2,a_3,b),\ y=f_2(a_1,a_2,a_3,b).\] The conditions that follow from possible connections between the arcs are: \[f_1(b,a_2,a_3,b)=f_2(b,a_2,a_3,b),\ \textrm{and}\ f_1(a_1,a_2,b,b)=f_2(a_1,a_2,b,b).\] An example of a pair of such functions is: 
\[f_1(a_1,a_2,a_3,b)=b\vee(a_1\wedge a_3)\ \textrm{and}\ f_2(a_1,a_2,a_3,b)=b\wedge(a_1\vee a_3).\]

\begin{figure}
\begin{center}
\includegraphics[height=3.2 cm]{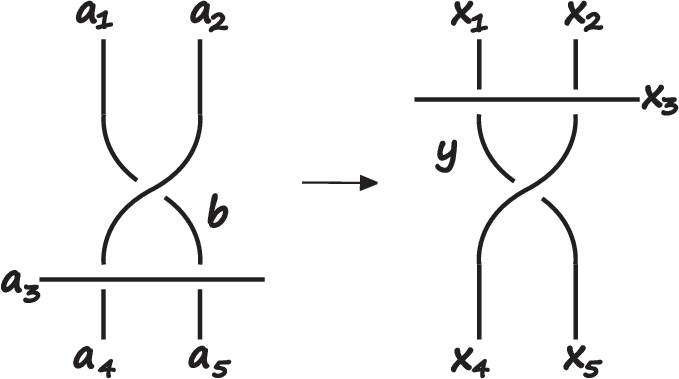}
\caption{}\label{lat3}
\end{center}
\end{figure}

Finally, consider the third labeled Reidemeister move R3, see Fig. \ref{lat3}. Here, 
\[x_i=f_i(a_1,\ldots,a_5,b),\ y=g(a_1,\ldots,a_5,b),\] where
\[f_i,\ g\colon \mathfrak{L}^6\to \mathfrak{L},\ \textrm{for}\ i=1,\ldots,5.\] 
The function $g$ could be arbitrary, and for the functions $f_i$, one needs the condition: if $a_i=a_j$, then
\[f_i(a_1,\ldots,a_i,\ldots,a_j,\ldots,a_5,b)=f_j(a_1,\ldots,a_i,\ldots,a_j,\ldots,a_5,b).\] 
Thus, for example, we could set $x_i=a_i$.

For decorated diagrams, the following version of the reachability problem is of an interest.

\begin{problem}[Weak reachability problem]
Let $\D$ be a category of diagrams and moves, with diagrams decorated by objects of a category $\E$.
Given a conditional category $\D\cap_f^w\C$, let $D_1^{E_1}\in(\D\cap_f^w\C)_0$, and let $D$ be the underlying diagram of one or more objects of $\D\cap_f^w\C$. Is there a finite sequence of elementary moves in $(\D\cap_f^w\C)_1$ leading from
$D_1^{E_1}$ to $D$? In other words: is there $D^{E}\in(\D\cap_f^w\C)_0$, for some $E$, such that $D^E$ is reachable from $D_1^{E_1}$? In such case we will write $D_1^{E_1}\rightarrow D$.
\end{problem}

\begin{problem} Let $\mathfrak{L}$ be a fixed lattice, and let $\D\cap_f^w\C$ be a conditional category of knot or link diagrams with arcs labeled by the elements of $\mathfrak{L}$.
Given two diagrams $D_1$ and $D_2$ with an incomplete label information, assign variables $x_1,x_2,\ldots,x_n$, 
to the arcs with no labels. The problem is to determine which elements of $\mathfrak{L}$, when substituted for the variables, and thus giving labelings $L_1$ and $L_2$, would make $D_2^{L_2}$ reachable from $D_1^{L_1}$. What are the properties of the set of solutions, and, in particular, how does it depend on the geometry of the diagrams, and on the choice of the lattice $\mathfrak{L}$? 

For a version of this question involving weak reachability, a partially labeled diagram $D_1$ would have variables assigned to it, and $D_2$ would be unlabeled, and with no variables. 
\end{problem}

\section{Knots with binary relations}
We will work with classical knot or link diagrams equipped with binary relations on the set of arcs, indicating which arcs can move over a given arc. We begin with some notation and terminology.

\begin{definition}
Let $S$ be a set, and let a binary relation $R\subseteq S\times S$ be given.
We will use the following notation:
\begin{itemize}
\item[--] $aR b$ or $R(a,b)=1$, to indicate that $(a,b)\in R$; 
\item[--] $a\cancel{R}b$ or $R(a,b)=0$ if $(a,b)\notin R$;
\item[--] $a \leadsto_{\ \ R} b$ if for any $c, d\in S$, $a R c$ implies $b R c$, and $d R a$ implies  $d R b$; then we say that $a$ {\it dominates} $b$ with respect to $R$. In particular, if $a\leadsto_{\ \ R} b$, then $a R a$ implies $b R b$;
\item[--] $a\sim_R b$ if $a\leadsto_{\ \ R} b$ and $b\leadsto_{\ \ R} a$. Then we say that $a$ and $b$ are {\it equivalent} with respect to $R$.
\end{itemize}
\end{definition}

\begin{remark}
For $S$ and $R$ as in the above definition, note that $\sim_R$ is an equivalence relation on $S$. Thus, we can form the quotient $A=S/\sim_R$, and introduce the induced relation on $A$, also denoted by $R$, by:
\[[a]R[b]\ \textrm{if and only if}\ aRb,\]
where $[x]$ denotes the equivalence class of $x$ with respect to $\sim_R$.
Also, we remark that the relation $\leadsto_{\ \ R}$ gives the preorder on $S$, and the poset on $S/\sim_R$, if we define:
\[[a]\leadsto_{\ \ R}[b]\ \textrm{if and only if}\ a\leadsto_{\ \ R}b.\]
\end{remark} 

Let $D$ be a classical knot or link diagram, $Arcs(D)$ denote the set of its arcs, and let $R$ be a binary relation on $Arcs(D)$; we will write $D^R$ to denote this situation. The question that immediately appears is: what should be the status of the new arcs created by the Reidemeister moves with respect to the relation $R$. There are many options, and they lead to different kinds of links; the choice can be influenced by desired applications. We choose to use Boolean terms, so that the resulting theory is broad, but still algorithmically manageable. We recall (closely following \cite{DP}) the definition of Boolean terms.
\begin{definition}
Let $A$ be the set of variables, and let $\vee$, $\wedge$, $\neg$, $0$, $1$ be the symbols used to axiomatize Boolean algebras. Then the class of Boolean terms, $\mathbf{BT}(A)$, is obtained as follows:
\begin{itemize}
\item[1)] $0$, $1\in\mathbf{BT}(A)$, and $a\in\mathbf{BT}(A)$ for all $a\in A$; 
\item[2)] if $p$, $q\in\mathbf{BT}(A)$, then $(p\vee q)$, $(p\wedge q)$, and $\neg p$ belong to $\mathbf{BT}(A)$;
\item[3)] every element of $\mathbf{BT}(A)$ is an expression formed by a finite number of applications of 1) and 2).
\end{itemize}
When elements of a given Boolean algebra $B$ are substituted for the variables of a Boolean term, $p(a_1,\ldots,a_n)$, an element of $B$ is obtained.
In particular, if $B$ is a two-element Boolean algebra $\mathbf{2}$ commonly used in logic, then every such $p$ defines a map $F_p\colon\mathbf{2}^n\to\mathbf{2}$. If $q(a_1,\ldots,a_n)$ is a Boolean term obtained from $p(a_1,\ldots,a_n)$ by the laws of Boolean algebra (in such case, we write $q\simeq p$), then $F_q=F_p$.
\end{definition}

\begin{lemma}\label{extendrel}
Let $A$ be a finite set with a binary relation $R$.
Then it is possible to extend $R$ onto $\mathbf{BT}(A)$, so that the following holds: for $p$, $p'$, $q$, $q'\in\mathbf{BT}(A)$, such that
$p\simeq p'$ and $q\simeq q'$, we have $p R q$ if and only if $p' R q'$. In other words, the relation between terms does not change if the laws of Boolean algebra are used to change the terms.
\end{lemma}
\begin{proof}
Suppose that $A$ has elements $a_1,\ldots,a_n$, and that $R(a_i,a_j)\in\{0,1\}$ is given for every $i$, $j\in\{1,\ldots,n\}$. The relation $R$ will be extended onto the set of Boolean terms $\mathbf{BT}(A)$ in the following steps. 

First, for all terms $p\in\mathbf{BT}(A)$, determine whether
\[a_i\ R\ p(a_1,a_2,\ldots,a_n)\] by taking the value $p(R(a_i,a_1),R(a_i,a_2),\ldots,R(a_i,a_n))\in\{0,1\}$. If $p'$ is such that
$p\simeq p'$, the truth tables for $p$ and $p'$ are the same. Thus, $a_i R p$ if and only if $a_i R p'$.

Next, for any $p$, $q\in\mathbf{BT}(A)$, decide whether \[q(a_1,a_2,\ldots,a_n)\ R\ p(a_1,a_2,\ldots,a_n)\]
by taking the value $q(R(a_1,p),R(a_2,p),\ldots,R(a_n,p))\in\{0,1\}$. If $q'\simeq q$, then $q R p$ if and only if $q' R p$, because we are looking at a particular row of two identical truth tables. Also, from the previous step, this row does not change if we replace $p$ by $p'$, such that $p\simeq p'$.

Finally, we note that the values $R(0,1)$ and $R(1,0)$ can be obtained as above, by writing $1=a_i\vee\neg a_i$ and $0=a_i\wedge\neg a_i$, for any $a_i$; we get $R(0,1)=0$ and $R(1,0)=1$. In the proof, we have made a choice of distributing $R$ first to the left, and then to the right, when determining if $qRp$, as is illustrated in the example that follows. 
\end{proof}

\begin{example}
Let the relation $R$ on the set of symbols $A=\{a,b,c\}$ be given by the following table.
\begin{center}
\begin{tabular}{c| c ccc} 
$R$&$a$ & $b$ & $c$ \\
\hline 
$a$ & 0 & 1 & 0 \\
$b$ & 1 & 0 & 1 \\
$c$ & 0 & 0 & 1 \\
\end{tabular}
\end{center}
Suppose that $p=(a\vee b)\wedge c$, and we want to determine if $pRp$, after the relation is extended onto the set of Boolean terms $\mathbf{BT}(A)$, as in the above proof.
\[(a\vee b)\wedge c\ R\ (a\vee b)\wedge c \Leftrightarrow \]
\[[(a\ R\ (a\vee b)\wedge c)\vee(b\ R\ (a\vee b)\wedge c)]\wedge (c\ R\ (a\vee b)\wedge c) \Leftrightarrow\]
\[[((aRa\vee aRb)\wedge aRc)\vee((bRa\vee bRb)\wedge bRc)]\wedge ((cRa\vee cRb)\wedge cRc) \Leftrightarrow\]
\[[((0\vee 1)\wedge 0)\vee((1\vee 0)\wedge 1)]\wedge ((0\vee 0)\wedge 1) \Leftrightarrow (0\vee 1)\wedge 0 \Leftrightarrow 0.\]
Thus, $p\cancel{R}p$.
\end{example}

\begin{remark}
Note that in Lemma \ref{extendrel}, the set $\mathbf{BT}(A)$ can be viewed as the free Boolean algebra freely generated by the elements of $A$. Then, the sublattice generated (using only $\vee$ and $\wedge$) by the same set, is the free distributive lattice, freely generated by the elements of $A$. For a given diagram $D^R$, denote the above Boolean algebra and distributive lattice generated by $Arcs(D)$, by $\mathcal{B}(D^R)$ and
$Distr(D^R)$, respectively.
\end{remark}

\begin{remark}
For $x\in\mathbf{BT}(A)$, and $F\subseteq\mathbf{BT}(A)$, let
\[x_{(R,F)}=\{y\in F\ | \ xRy\}\]
\[x^{(R,F)}=\{y\in F\ | \ yRx\}.\]
Then, if $p$ is any Boolean term with the set of variables $A=\{a_1,\ldots,a_n\}$, the following holds:
\[p(a_1,\ldots,a_n)_{(R,F)}=p((a_1)_{(R,F)},\ldots,(a_n)_{(R,F)}),\ \textrm{for any}\ F,\]
\[p(a_1,\ldots,a_n)^{(R,A)}=p((a_1)^{(R,A)},\ldots,(a_n)^{(R,A)}),\]
where on the right hand side there are expressions obtained by replacing $\vee$, $\wedge$, and $\neg$, by the union, intersection, and complement of sets, respectively.
For example:
\[((a_1\vee a_3)\wedge a_6)_{(R,F)}=((a_1)_{(R,F)}\cup(a_3)_{(R,F)})\cap(a_6)_{(R,F)}.\]
\end{remark}

\begin{lemma}\label{dominates}
Suppose $(A,R)$ is a set with relation, and $Distr(A^R)$ is the free distributive lattice on $A$ with extended relation obtained as in Lemma \ref{extendrel}.
If $p$, $q\in Distr(A^R)$ and $p\leq q$ in the lattice order, i.e. $p\wedge q=p$, then $p\leadsto_{\ \ R} q$.
\end{lemma}
\begin{proof}
We need to show that: (i) for any $c\in (Distr(A),R)$, $p R c$ implies $q R c$, and (ii) for any $d\in (Distr(A),R)$, 
$d R p$ implies $d R q$. For (i), we note that: 
\[p R c \ \textrm{iff}\ p\wedge q R c \ \textrm{iff}\  p R c \wedge q R c,\]
and the last logical expression implies $q R c$.

For (ii), let $d=d(a_1,\ldots,a_n)$. Then we have:
\[R(d(a_1,\dots,a_n),p)=1 \ \textrm{iff}\ d(R(a_1,p),\ldots, R(a_n,p))=1 \ \textrm{iff}\]
\[d(R(a_1,p\wedge q),\ldots,R(a_n,p\wedge q))=1 \ \textrm{iff}\]
\[d(R(a_1,p)\wedge R(a_1,q),\ldots,R(a_n,p)\wedge R(a_n,q))=1.\] 
It is a standard fact that lattice terms are isotone, i.e., if $s(x_1,\ldots,x_n)$ is a lattice term with variables $x_1,\ldots,x_n$, and
$b_1\leq c_1,\ldots, b_n\leq c_n$, then $s(b_1,\ldots,b_n)\leq s(c_1,\ldots,c_n).$
Here, we have:
\[R(a_1,p)\wedge R(a_1,q)\leq R(a_1,q), \ldots, R(a_n,p)\wedge R(a_n,q)\leq R(a_n,q),\]
and thus,
\[1=d(R(a_1,p)\wedge R(a_1,q),\ldots,R(a_n,p)\wedge R(a_n,q))\leq d(R(a_1,q),\ldots,R(a_n,q)),\]
so $d(R(a_1,q),\ldots,R(a_n,q))=1$, that is, $R(d,q)=1$.
\end{proof}

\begin{figure}
\begin{center}
\includegraphics[height=10 cm]{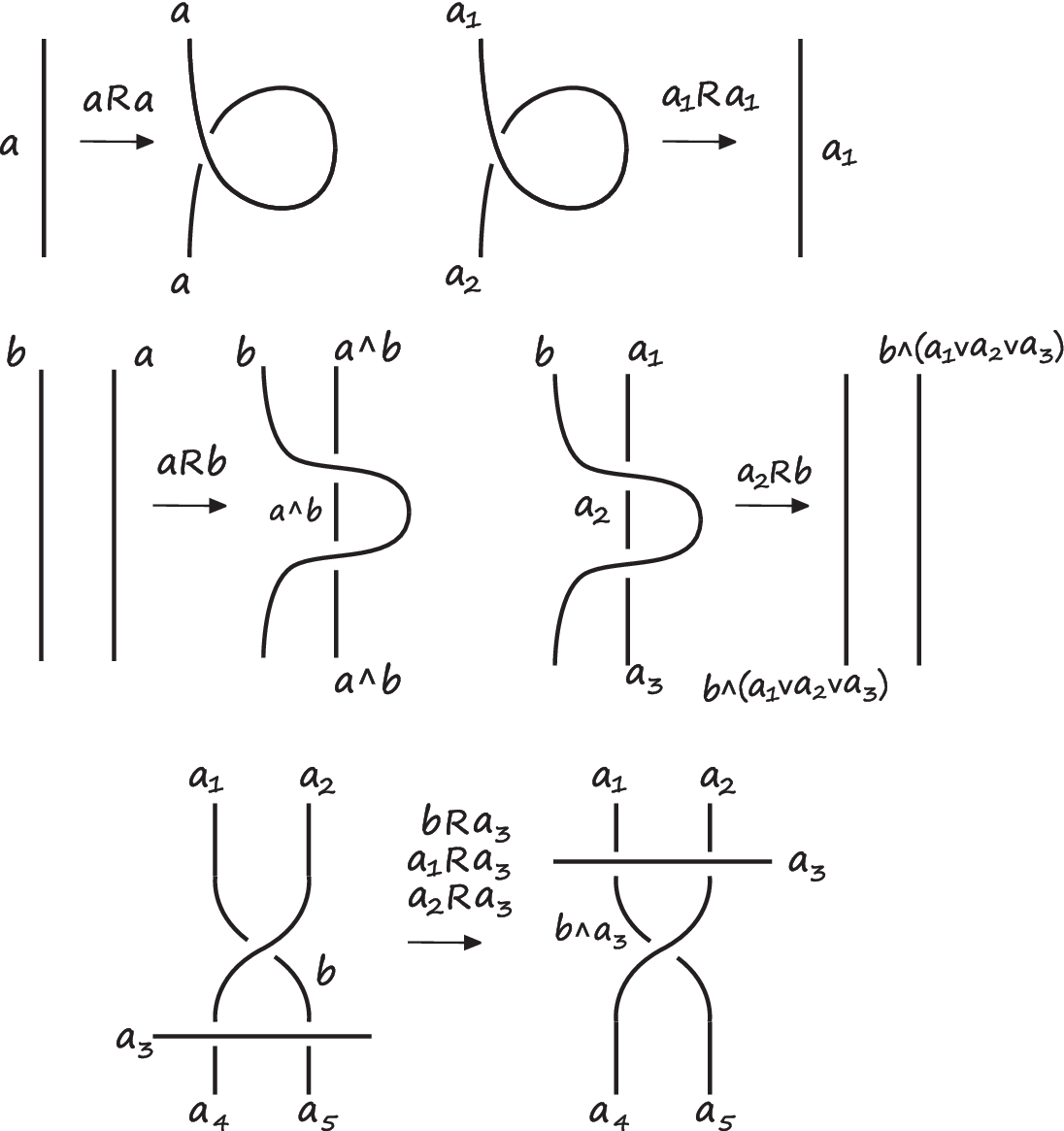}
\caption{}\label{moveset}
\end{center}
\end{figure}

Given these preliminaries, we are ready to define the category of links with relations.

\begin{definition}
We will define a relation on arcs of the diagram obtained after a Reidemeister move. Let $D_1^{R_1}$ denote a knot or link diagram $D_1$ with a binary relation $R_1$ on $Arcs(D_1)$. Let $D_2$ be a diagram obtained by performing a Reidemeister move on $D_1$. 
The relation $R_2$ on $Arcs(D_2)$ is defined by assigning elements of $\mathcal{B}(D_1^{R_1})$ to the arcs of 
$D_2$. We assume that if an arc $a$ is not involved in the move (and thus can be viewed as being in the intersection
$Arcs(D_1)\cap Arcs(D_2)$), then the symbol assigned to it is $a$. If $b_1[t_1]$ and $b_2[t_2]$
denote the arcs $b_1$, $b_2\in Arcs(D_2)$ with the terms $t_1$, $t_2\in\mathcal{B}(D_1^{R_1})$ assigned to them, then we set:
\[b_1 R_2 b_2\ \textrm{iff}\ t_1 R_1 t_2.\]
To have consistency in defining relations after the moves, we label the arcs involved in the Reidemeister moves using Boolean terms, to indicate how the relation is changing when going from $D_1$ to $D_2$. We call such labeled moves {\it relation Reidemeister moves}. One such set of moves is illustrated in Figure \ref{moveset}. 
It shows the labels $t\in\mathcal{B}(D_1)$ assigned to the arcs of $D_2$.
When a set of relation Reidemeister moves is chosen, we obtain a category of diagrams and moves, denoted by 
$\D^{\mathfrak{R}}$.
\end{definition}

\begin{remark}
We note that, in general, one could take a set of moves containing several moves of each Reidemeister type, differing by the relation (or labeling) that they induce.
\end{remark}
 
\begin{remark}
We defined categories $\D^{\mathfrak{R}}$ in order to introduce conditional categories with the restriction depending on the relation graph of $R$ in a given $D^R$. Here we analyze one of the simplest conditions: 
\begin{equation}\label{basiccond}
\textrm{for}\ a_1,\, a_2\in Arcs(D^R),\ a_2\ \textrm{can move over}\ a_1\ \textrm{if}\ a_1 R\, a_2. \tag{C1}
\end{equation}

We mention that the categories with this condition can be expressed as in Definition \ref{genknot}. Take $\C$ to be the category of sets with inclusions (i.e., $\exists m\colon S_1\to S_2$ iff $S_1\subseteq S_2$). If $D^R$ is a diagram before the move, let $c\in Arcs(D^R)$ be the arc that is going to move over arcs, say $a_1\in Arcs(D^R)$ if the move is of type I or II, and $a_1$, $a_2$, $a_3\in Arcs(D^R)$, in the case of the third Reidemeister move. Let 
\[A_i=\{b\in Arcs(D^R)\,|\, a_i R b\},\]
for $i\in\{1,2,3\}$. In the case of the first or second Reidemeister move take:
\[f(D_1,D_2)=\{(\{c\},A_1),(\{c\},A_1),(\{c\},A_1)\},\] and if the third Reidemeister move is performed, let
\[f(D_1,D_2)=\{(\{c\},A_1),(\{c\},A_2),(\{c\},A_3)\}.\]
Then the formula 
\[w=(\exists m_1\colon Obj_1\to Obj_1')\wedge(\exists m_2\colon Obj_2\to Obj_2')\wedge(\exists m_3\colon Obj_3\to Obj_3')\] allows to write these conditional categories as $\D^{\mathfrak{R}}\cap_f^w\C$. 
\end{remark} 

The conditions above arrows in Fig. \ref{moveset} are used after passing to the conditional theory 
$\D^{\mathfrak{R}}\cap_f^w\C$ (in which the Condition (\ref{basiccond}) holds). Then, some of the moves in Fig. \ref{moveset} might be irreversible (depending on the relation in $D^R$), as the following example illustrates.

\begin{figure}
\begin{center}
\includegraphics[height=14 cm]{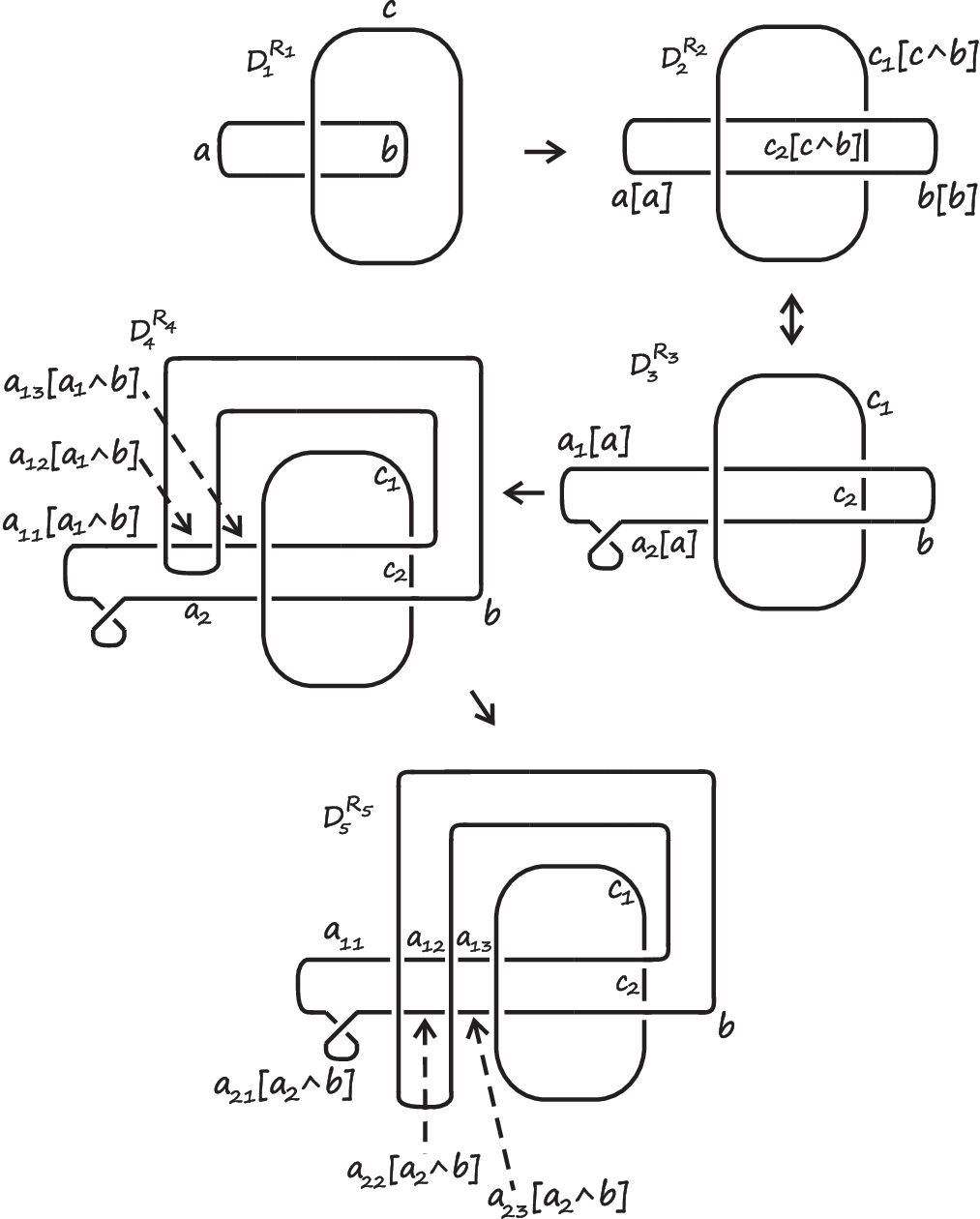}
\caption{}\label{reldiagrams}
\end{center}
\end{figure}

\begin{example}
Figure \ref{reldiagrams} shows an example of relation Reidemeister moves from Fig. \ref{moveset} applied to diagrams with relations. The following tables describe the initial and subsequent relations assigned to the diagrams. The relation table of $R_5$ (not shown here) contains only zeros, which means that no further Reidemeister moves are possible; $D_5^{R_5}$ is a terminal state. 

\begin{center}
\begin{tabular}{c| c ccc} 
$R_1$ & $a$ & $b$ & $c$ \\
\hline 
$a$ & 1 & 1 & 0 \\
$b$ & 0 & 0 & 0 \\
$c$ & 1 & 1 & 1 \\
\end{tabular}
\end{center}

\begin{center}
\begin{tabular}{c| c cccc} 
$R_2$ &$a$ & $b$ & $c_1$ & $c_2$\\
\hline 
$a$   & 1 & 1 & 0 & 0 \\
$b$   & 0 & 0 & 0 & 0 \\
$c_1$ & 0 & 0 & 0 & 0 \\
$c_2$ & 0 & 0 & 0 & 0 \\
\end{tabular}
\end{center}

\begin{center}
\begin{tabular}{c| c ccccc} 
$R_3$ & $a_1$ & $a_2$ & $b$ & $c_1$ & $c_2$\\
\hline 
$a_1$   & 1 & 1 & 1 & 0 & 0\\
$a_2$   & 1 & 1 & 1 & 0 & 0\\
$b$   & 0 & 0 & 0 & 0 & 0\\
$c_1$ & 0 & 0 & 0 & 0 & 0\\
$c_2$ & 0 & 0 & 0 & 0 & 0\\
\end{tabular}
\end{center}

\begin{center}
\begin{tabular}{c| c ccccccc} 
$R_4$ & $a_{11}$ & $a_{12}$ & $a_{13}$ & $a_2$ & $b$ & $c_1$ & $c_2$\\
\hline 
$a_{11}$ & 0 & 0 & 0 & 0 & 0 & 0 & 0\\
$a_{12}$ & 0 & 0 & 0 & 0 & 0 & 0 & 0\\
$a_{13}$ & 0 & 0 & 0 & 0 & 0 & 0 & 0\\
$a_2   $ & 1 & 1 & 1 & 1 & 1 & 0 & 0\\
$b     $ & 0 & 0 & 0 & 0 & 0 & 0 & 0\\
$c_1   $ & 0 & 0 & 0 & 0 & 0 & 0 & 0\\
$c_2   $ & 0 & 0 & 0 & 0 & 0 & 0 & 0\\
\end{tabular}
\end{center}

\end{example}

\begin{lemma}\label{biggerrel}
Let $R_1\subseteq R_2$, and $D^{R_1}\rightarrow D_3^{R_3}$, in a given $\D^{\mathfrak{R}}\cap_f^w\C$ possessing a set of moves that use only labels involving $\wedge$ and $\vee$.
Then $D^{R_2}\rightarrow D_3^{R_4}$, where $R_3\subseteq R_4$.
\end{lemma}
\begin{proof} We compare the extended relations on the distributive lattices:\\ 
$Distr(D^{R_1})$ and $Distr(D^{R_2})$.
We prove that for any terms $p=p(a_1,\ldots,a_n)$ and $q=q(a_1,\ldots,a_n)$, $p R_1 q$ implies $p R_2 q$.
Indeed, suppose 
\[R_1(p(a_1,\ldots,a_n),q(a_1,\ldots,a_n))=1;\] it is equivalent to the following:
\[p(R_1(a_1,q(a_1,\ldots,a_n)),\ldots,R_1(a_n, q(a_1,\ldots,a_n)))=1\ \textrm{iff}\]
\[p(q(R_1(a_1,a_1),\ldots,R_1(a_1,a_n)),\ldots,q(R_1(a_n,a_1),\ldots,R_1(a_n,a_n)))=1.\]
$R_1(a_i,a_j)\leq R_2(a_i,a_j)$, for any generators $a_i$ and $a_j$. Thus, from isotonity of $q$, we have:
\[q(R_1(a_1,a_1),\ldots,R_1(a_1,a_n))\leq q(R_2(a_1,a_1),\ldots,R_2(a_1,a_n)),\ldots\]
\[q(R_1(a_n,a_1),\ldots,R_1(a_n,a_n))\leq q(R_2(a_n,a_1),\ldots,R_2(a_n,a_n)).\]
Therefore, from the isotonity of $p$:
\[1=p(q(R_1(a_1,a_1),\ldots,R_1(a_1,a_n)),\ldots,q(R_1(a_n,a_1),\ldots,R_1(a_n,a_n)))\leq\]
\[p(q(R_2(a_1,a_1),\ldots,R_2(a_1,a_n)),\ldots,q(R_2(a_n,a_1),\ldots,R_2(a_n,a_n)))\ \textrm{,i.e.,}\] 
\[p(R_2(a_1,q),\ldots,R_2(a_n,q))=1 \ \textrm{,i.e.,}\ R_2(p,q)=1.\]
It follows that the relation $R$ on the arcs of a diagram obtained from $D^{R_1}$ via a Reidemeister move, is included in the relation $R'$ on the arcs of a diagram obtained from $D^{R_2}$ via the same move. Thus, any move permitted by $R$ is also permitted by $R'$. The situation repeats itself after every Reidemeister move, including the step when $D_3$ is reached.
\end{proof}

\begin{definition}
We say that a set of moves of $\D^{\mathfrak{R}}$ is {\it entropy decreasing} if for every $t\in\mathcal{B}(D_1)$ that is used to label an arc of $D_2$ after a Reidemeister move,
there is an arc $a(t)\in Arcs(D_1)$ (viewed as an element of $\mathcal{B}(D_1)$) such that $t\leadsto_{\ \ {R_1}} a(t)$.
\end{definition}

It follows from Lemma \ref{dominates} that the set of moves in Fig \ref{moveset} is entropy decreasing.

\begin{definition}
Let $A^R$ and $X^\R$ be two sets with relations. A monotone map is a map $f\colon A^R\to X^\R$ such that
if $a_1 R a_2$ then $f(a_1)\R f(a_2)$, for any $a_1$, $a_2\in A$.
\end{definition}

\begin{theorem} \label{entrdec}
Let $\D^{\mathfrak{R}}\cap_f^w\C$ has an entropy decreasing set of moves. Let $X^{\R}$ be a set with relation. Let $D_1^{R_1}$ and $D_2^{R_2}$ be such that $D_1^{R_1}\rightarrow D_2^{R_2}$.
Then for every monotone map $f_1\colon Arcs(D_1)^{R_1}\to X^{\R}$, there is a monotone map
$f_2\colon Arcs(D_2)^{R_2}\to X^{\R}$ such that $Im(f_2)\subseteq Im(f_1)$, where $Im(f)$ denotes the image of a map $f$.
\end{theorem}
\begin{proof}
Suppose that $D_2$ is obtained from $D_1$ by a single Reidemeister move.
Given a monotone map $f_1\colon Arcs(D_1)^{R_1}\to X^{\R}$, define
\[f_2\colon Arcs(D_2)^{R_2}\to X^{\R}\ \textrm{by}\ f_2(d)=f_1(a(t)),\]
where $t\in \mathcal{B}(D_1)$ is the label assigned to $d\in Arcs(D_2)$, and $a(t)\in Arcs(D_1)$ is such that
$t\leadsto_{\ \ {R_1}} a(t)$.

Let $d_1$, $d_2\in Arcs(D_2)$, and let $t_1$, $t_2\in\mathcal{B}(D_1)$ be the labels assigned to them
after a Reidemeister move. Recall that $d_1 R_2 d_2 \ \textrm{iff}\ t_1 R_1 t_2$.
We need to check that $f_2$ is monotone with respect to $R_2$. Suppose $d_1 R_2 d_2$, and so $t_1 R_1 t_2$.
The map $f_2$ is monotone if $f_2(d_1)\R f_2(d_2)$, i.e., $f_1(a(t_1))\R f_1(a(t_2))$.
Since $t_1\leadsto_{\ \ {R_1}} a(t_1)$, we have:
\[t_1 R_1 t_2 \ \textrm{implies}\ a(t_1)R_1 t_2.\]
Since $t_2\leadsto_{\ \ {R_1}} a(t_2)$,
\[a(t_1) R_1 t_2 \ \textrm{implies}\ a(t_1)R_1 a(t_2).\]
Then, since $f_1$ is $R_1$-monotone, we have $f_1(a(t_1))\R f_1(a(t_2))$, proving the $R_2$-monotonicity of $f_2$.
From the definition of $f_2$ follows that $Im(f_2)\subseteq Im(f_1)$. The above proof was for $D_1^{R_1}$ and $D_2^{R_2}$ differing by a single Reidemeister move, but we note that for a longer sequence of moves, a monotone map after each move can be constructed as above from a monotone map before the move, and the condition regarding the images is satisfied.
\end{proof}

Now we take a look at a behavior of the standard properties of relations under the moves, and under certain extensions.

\begin{lemma}\label{varprops}
Let $R$ be a relation on a set $X$, and take its extension, also denoted by $R$, to $\mathcal{B}(X^R)$, as in Lemma \ref{extendrel}.
\begin{enumerate}
\item[(i)] If $R$ is reflexive on $X$, then $R$ is also reflexive on any set $X\cup S$, where $S$ contains elements
of $\mathcal{B}(X^R)$ that can be written using only the operation $\vee$, or $S$ has only elements that are negations of elements of $X$;
\item[(ii)] If $R$ is transitive on $X$, then $R$ is also transitive on any set $X\cup S$, where $S$ contains elements
of $\mathcal{B}(X^R)$ that can be written using only the operation $\wedge$;
\item[(iii)] If $R$ is symmetric on $X$, then $R$ is also symmetric on a set $X\cup \{t\}$, where $t$ is any element
of $\mathcal{B}(X^R)$;
\item[(iv)] If $R$ is symmetric on $X$, then $R$ is also symmetric on any set $X\cup S$, where $S$ contains elements
of $\mathcal{B}(X^R)$ that can be written using only the operation $\wedge$, or $S$ has only elements that can be written using only $\vee$.
\end{enumerate}
\end{lemma}
\begin{proof}
(i) If $t=\bigvee_{i=1}^k a_i$, $a_i\in X$, then
\[\bigvee_{i=1}^k a_i\,R\,\bigvee_{i=1}^k a_i\ \textrm{iff}\ \bigvee_{i,j\in\{1,\ldots,k\}} a_i R a_j.\]
Since $a_i R a_i$, the above disjunction is true, that is, $t R t$.

Regarding the negations, we note that:
\[\neg a R \neg a\ \textrm{iff}\ \neg(a R \neg a)\ \textrm{iff}\ \neg\neg(a R a)\ \textrm{iff}\ aRa.\]

\noindent (ii) If $p=\bigwedge_{i=1}^{n}a_i$, $q=\bigwedge_{j=1}^{m}b_j$, and $s=\bigwedge_{k=1}^{l}c_k$,
with all $a_i$, $b_j$, $c_k\in X$; then:
\[\bigwedge_{i=1}^{n}a_i\,R\,\bigwedge_{k=1}^{l}c_k\ \textrm{iff}\ 
\bigwedge_{i\in\{1,\ldots,n\},k\in\{1,\ldots,l\}} a_i R c_k.\]
We have:
\[p R q\ \textrm{iff}\ \bigwedge_{i\in\{1,\ldots,n\},j\in\{1,\ldots,m\}} a_i R b_j,\ \textrm{and}\]
\[q R s\ \textrm{iff}\ \bigwedge_{j\in\{1,\ldots,m\},k\in\{1,\ldots,l\}} b_j R c_k,\]
thus, from transitivity of $R$ on generators, $\bigwedge_{i\in\{1,\ldots,n\},k\in\{1,\ldots,l\}} a_i R c_k$ holds, 
that is, $p R s$.

\noindent (iii) Suppose that $t=t(a_1,\ldots,a_n)$ is a term from $\mathcal{B}(X^R)\setminus X$, and $b\in X$.
Then we have:
\[t(a_1,\ldots,a_n)\,R\,b \ \textrm{iff}\ R(t(a_1,\ldots,a_n),b)=1\ \textrm{iff}\ 
t(R(a_1,b),\ldots,R(a_n,b))=1\ \textrm{iff}\]
\[t(R(b,a_1),\ldots,R(b,a_n))=1\ \textrm{iff}\ R(b,t(a_1,\ldots,a_n))=1\ \textrm{iff}\ b\,R\,t(a_1,\ldots,a_n).\]

\noindent (iv) Let $p=\bigvee_{i=1}^{n}a_i$, $q=\bigvee_{j=1}^{m}b_j$, where
all $a_i$, $b_j\in X$, then:
\[\bigvee_{i=1}^{n}a_i\,R\,\bigvee_{j=1}^{m}b_j\ \textrm{iff}\ 
\bigvee_{i\in\{1,\ldots,n\},j\in\{1,\ldots,m\}} a_i R b_j\ \textrm{iff}\]
\[\bigvee_{i\in\{1,\ldots,n\},j\in\{1,\ldots,m\}} b_j R a_i\ \textrm{iff}\ 
\bigvee_{j=1}^{m}b_j\,R\,\bigvee_{i=1}^{n}a_i.\]
The proof for $\wedge$ is analogous.
\end{proof}

\begin{theorem}\label{preserve}
Consider a conditional theory $\D^{\mathfrak{R}}\cap_f^w\C$ that has the set of moves illustrated in 
Fig. \ref{moveset}. Suppose that $D_1^{R_1}\rightarrow D_2^{R_2}$, and that the relation $R_1$ is symmetric on $Arcs(D_1)$.
Then the following holds:
\begin{enumerate}
\item[(i)] $R_2$ is symmetric on $Arcs(D_2)$;
\item[(ii)] If $R_1$ is reflexive on $Arcs(D_1)$, then $R_2$ is reflexive on $Arcs(D_2)$;
\item[(iii)] If $R_1$ is transitive on $Arcs(D_1)$, then $R_2$ is transitive on $Arcs(D_2)$.
\end{enumerate}
\end{theorem} 
\begin{proof}
Part (i) follows from Lemma \ref{varprops}(iii), because each move uses at most one label from $Distr(D_1^{R_1})\setminus Arcs(D_1)$, so the relation after a Reidemeister move, and subsequent relations, remain symmetric.

Parts (ii) and (iii) use conditionality. To prove (ii) we need to check that if 
$a$, $a_1$, $a_2$, $a_3$, $b\in Arcs(D_1)$, and $a R_1 b$, $a_2 R_1 b$, then $a\wedge b \,R_1\, a\wedge b$, and
$b\wedge (a_1\vee a_2\vee a_3) \,R_1\, b\wedge (a_1\vee a_2\vee a_3)$.
We have:
\[R_1(a\wedge b,a\wedge b)=1\ \textrm{iff}\ R_1(a,a\wedge b)=1\wedge R_1(b,a\wedge b)=1\ \textrm{iff}\]
\[R_1(a,a)=1 \wedge R_1(a,b)=1 \wedge R_1(b,a)=1 \wedge R_1(b,b)=1.\]
The first and fourth parts of the above conjunction are true because of the reflexivity of $R_1$ on $Arcs(D_1)$, the second and third follow from conditionality and symmetry.

Let us investigate the second element:
\[b\wedge (a_1\vee a_2\vee a_3) \,R_1\, b\wedge (a_1\vee a_2\vee a_3)\ \textrm{iff}\]
\[[b\,R_1\, b\wedge (a_1\vee a_2\vee a_3)] \wedge [a_1\vee a_2\vee a_3 \,R_1\, b\wedge (a_1\vee a_2\vee a_3)]
\ \textrm{iff}\ \] 
\[[(b\,R_1\, b) \wedge (b\,R_1\,(a_1\vee a_2\vee a_3))]\wedge[(a_1\,R_1\, b\wedge (a_1\vee a_2\vee a_3))\]
\[\vee (a_2\,R_1\, b\wedge (a_1\vee a_2\vee a_3)) \vee (a_3\,R_1\, b\wedge (a_1\vee a_2\vee a_3))]\ \textrm{iff}\]
\[[(b\,R_1\, b) \wedge (b\,R_1\,a_1 \vee b\,R_1\,a_2 \vee b\,R_1\,a_3)]\wedge 
[((a_1\,R_1\, b)\wedge (a_1\,R_1\, a_1\vee a_2\vee a_3))\] 
\[\vee ((a_2\,R_1\, b)\wedge (a_2\,R_1\, a_1\vee a_2\vee a_3)) \vee ((a_3\,R_1\, b)\wedge (a_3\,R_1\, a_1\vee a_2\vee a_3))].\]
The true value of the last logical sentence follows from $R_1(b,b)=1$, $R_1(b,a_2)=1$ (conditionality and symmetry), and the truth of the part \[(a_2\,R_1\, b)\wedge (a_2\,R_1\, a_1\vee a_2\vee a_3),\] which follows from conditionality and reflexivity. 

In the proof of part (iii), there are some cases. First, we note that because of Lemma \ref{varprops}(ii), adding $a\wedge b$ to the set of labels preserves the transitivity. Now, suppose that 
\[c\,R_1\, b\wedge (a_1\vee a_2\vee a_3)\ \textrm{and}\ b\wedge (a_1\vee a_2\vee a_3)\,R_1\,d,\]
where $c$, $d\in Arcs(D_1)$. Then it follows that $c R_1 b$ and $b R_1 d$, which gives $c R d$ from the transitivity of
$R_1$ on $Arcs(D_1)$. If
\[b\wedge (a_1\vee a_2\vee a_3)\,R_1\,c \ \textrm{and}\ c\,R_1\,d,\ \textrm{then}\ \]
\[(b\,R_1\,c) \wedge (a_1\vee a_2\vee a_3\,R_1\,c),\ \textrm{i.e.,}\ \]
\[(b\,R_1\,c) \wedge (a_1\,R_1\,c \vee a_2\,R_1\,c \vee a_3\,R_1\,c).\] 
Since $a_2 R_1 b$, $b R_1 c$, and $R_1$ is transitive on $Arcs(D_1)$, we have $a_2 R_1 c$. Thus, from transitivity,
$b R_1 d$ and $a_2 R_1 d$. Therefore,
\[(b\,R_1\,d) \wedge (a_1\vee a_2\vee a_3\,R_1\,d),\ \textrm{i.e.,}\ \]
\[b\wedge (a_1\vee a_2\vee a_3)\,R_1\,d.\]
The proof of transitivity when
\[c\,R_1\,d \ \textrm{and}\ d\,R_1\,b\wedge (a_1\vee a_2\vee a_3)\]
is very similar, although it additionally uses symmetry. Finally, we consider:
\[b\wedge (a_1\vee a_2\vee a_3)\,R_1\,c \ \textrm{and}\ c\,R_1\,b\wedge (a_1\vee a_2\vee a_3)\]
The reflexivity condition:
\[b\wedge (a_1\vee a_2\vee a_3)\,R_1\,b\wedge (a_1\vee a_2\vee a_3),\]
that we need here, is expanded in the proof of (ii). It holds, because $b R_1 c$ and $c R_1 b$ implies $b R_1 b$;
also $a_2 R_1 b$ and $b R_1 a_2$ yields $a_2 R a_2$.
\end{proof}

\begin{definition}
We will say that a relation $R$ is $\it semi-transitive$ if $a R b$ and $b R c$ implies $a R c$, for any distinct
elements $a$, $b$, $c$. 
\end{definition}
Taking a semi-transitive closure, instead of the usual transitive closure, avoids the creation of reflexivity in the presence of symmetry.

Now we can combine Lemma \ref{biggerrel}, Theorem \ref{entrdec}, and Theorem \ref{preserve}, to obtain the following
corollary.

\begin{corollary}
Let $\D^{\mathfrak{R}}\cap_f^w\C$ be the conditional theory with the set of moves illustrated in 
Fig. \ref{moveset}. For a given relation $R$, let $R^{ST}$ denote the smallest symmetric and semi-transitive relation containing $R$.
Let $D_1^{R_1}$ and $D_2^{R_2}$ be such that $D_1^{R_1}\rightarrow D_2^{R_2}$. Let $X^{\R}$ be a set with relation.
Then, for every monotone map $f_1\colon Arcs(D_1)^{R_1^{ST}}\to X^{\R}$, there is a monotone map
$f_2\colon Arcs(D_2)^{R_2^{ST}}\to X^{\R}$ such that $Im(f_2)\subseteq Im(f_1)$.
\end{corollary}
\begin{proof}
Since $R_1\subseteq R_1^{ST}$, from Lemma \ref{biggerrel} there exists a relation $R_3$ such that
$D_1^{R_1^{ST}}\rightarrow D_2^{R_3}$, with $R_2\subseteq R_3$. Theorem \ref{preserve} implies that the symmetry and semi-transitivity are preserved by the moves, thus, $R_3$ is a symmetric and semi-transitive relation on $Arcs(D_2)$, and has to contain $R_2^{ST}$.
Let $f_1\colon Arcs(D_1)^{R_1^{ST}}\to X^{\R}$ be monotone. As we already mentioned, the set of moves in Fig. \ref{moveset}
is entropy decreasing. Then, from Theorem \ref{entrdec} there exists $f_2\colon Arcs(D_2)^{R_3}\to X^{\R}$ that is monotone and such that $Im(f_2)\subseteq Im(f_1)$. Then, since $R_2^{ST}\subseteq R_3$, the same $f_2$ is monotone on 
$Arcs(D_2)^{R_2^{ST}}$.
\end{proof}

\section{Links with relations on the set of components}

\begin{figure}
\begin{center}
\includegraphics[height=6.6cm]{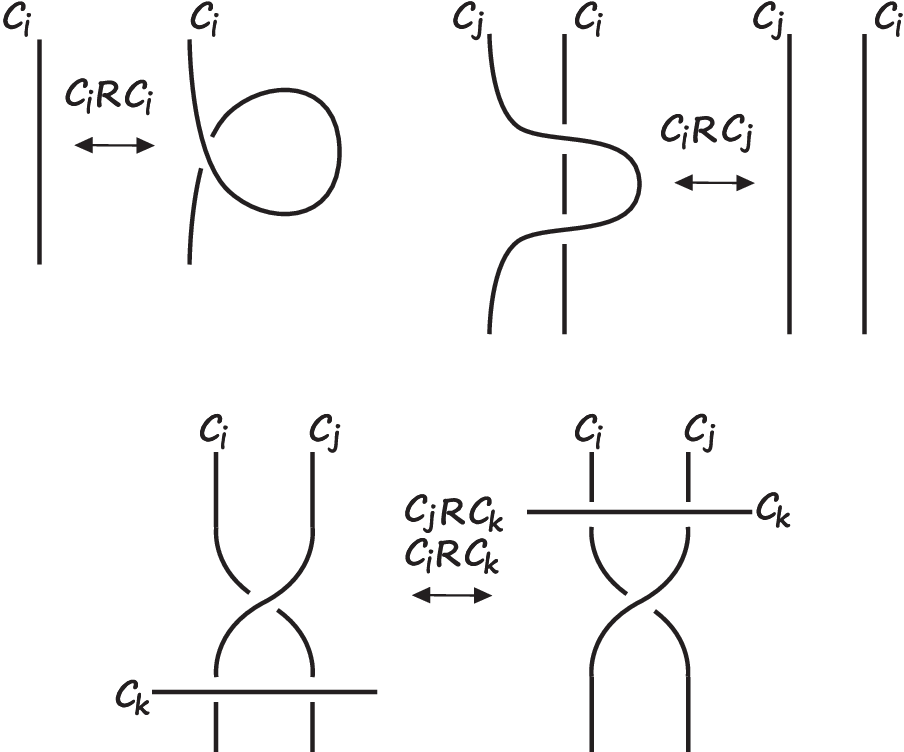}
\caption{}\label{compmoves}
\end{center}
\end{figure}

Let $L$ be an oriented link diagram with the set of components $C=\{C_1,\ldots,C_s\}$, and let $R$ be a binary relation on $C$.
We impose the condition: the component $C_j$ can move over the component $C_i$ if and only if $C_i\,R\,C_j$, for $i$, 
$j\in\{1,\ldots,s\}$. In terms of the Reidemeister moves, this condition is depicted in Fig. \ref{compmoves}, with the conditions required for the move indicated above the arrows. The moves, if allowed, are reversible, and thus we will be looking for invariants of links with relations.

\begin{definition}
Suppose that the component $C_i$ has arcs $a_1$, $a_2,\ldots,a_k$, and the component $C_j$ has arcs $b_1$, $b_2,\ldots,b_n$. Then we define the induced relation on arcs by:
\[a_p\,R\,b_q\ \textrm{iff}\ C_i\,R\,C_j\ \textrm{for}\ p\in\{1,\ldots,k\}, q\in\{1,\ldots,n\},
\ i,j\in\{1,\ldots,s\}.\]
In particular, $a_p\,\sim_R\,a_{p'}$, for $p$, $p'\in\{1,\ldots,k\}$, and $b_q\,\sim_R\,b_{q'}$, for $q$, $q'\in\{1,\ldots,n\}$. In other words, the arcs belonging to the same component are equivalent with respect to the relation $R$. Also, we assume that after a Reidemeister move, the relation on the new set of arcs is again induced from the relation on the set of components.
\end{definition}

The structure that will be of use to us is closely related to racks and quandles. For completeness, we will give the definition and examples of racks and quandles; see, for example, \cite{CJKLS, FRS, Joy} for more about these structures.

\begin{definition} \label{quandle}
A {\it quandle}, $X$, is a set with two binary operations: 
$(x,y) \mapsto x*y$, and $(x,y) \mapsto x\,\bar{*}\,y$ such that:
\begin{enumerate}
\item[(1)] $x*x=x=x\,\bar{*}\,x$, for any $x\in X$,
\item[(2)] $x*y\,\bar{*}\,y=x=x\,\bar{*}\,y*y$, for any $x$, $y\in X$, 
\item[(3)] $(x*y)*z=(x*z)*(y*z)$, for any $x$, $y$, $z\in X$. 
\end{enumerate}
\end{definition}
\begin{definition}
A {\it rack} is a structure satisfying the axioms (2) and (3) of Def. \ref{quandle}.
\end{definition}
In general, racks and quandles are nonassociative structures, and the operations are performed from left to right if their order is not indicated by brackets.
 
The following examples often appear in applications to knot theory:

\begin{enumerate}
\item[-]
Any group $G$ is a quandle with operations: 
\[a*b=b^{-n} a b^n,\ a\,\bar{*}\,b=b^{n} a b^{-n},\] where $n$ is any integer. 
\item[-]
Another quandle is obtained from a given group $G$ if 
\[a*b=a\,\bar{*}\,b=b a^{-1} b.\] 
It is called the core quandle of $G$. 
\item[-]
Let $X=\{0,1,\ldots,n-1\}$. Then, a quandle, called the dihedral quandle, is obtained by using operations:
\[i*j=i\,\bar{*}\,j=2j-i \pmod{n},\ \textrm{for}\ i, j\in X.\]
\item[-] An example of a rack that in general is not a quandle is a $G$-set $X$ with operations:
\[x*y=x\cdot g,\ \textrm{and}\ x\,\bar{*}\,y=x\cdot g^{-1},\]
where, $x$, $y\in X$, and $g$ is a fixed element of a group $G$.
\end{enumerate}

In classical knot theory there is the following notion of a quandle coloring of a link diagram.

\begin{figure}
\begin{center}
\includegraphics[height=1.5 cm]{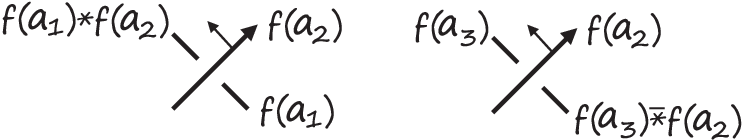}
\caption{The quandle coloring rule for knot and link diagrams.}\label{qcols}
\end{center}
\end{figure}

\begin{definition}\label{colorings}
Let $X$ be a fixed quandle, and let $Arcs(D)$ be the set of arcs of the link diagram $D$. The normals to the arcs are given in such a way that the pair (tangent vector, normal vector) matches the usual orientation of the plane. A {\it quandle coloring} is a map $f\colon Arcs(D)\to X$ such that at every crossing, the relation depicted 
in Fig. \ref{qcols} is satisfied. More precisely, if $a_2$ is an over-arc at a crossing, and $a_1$, $a_3$ are the under-arcs such that the normal of the over-arc points from $a_1$ to $a_3$, then it is required that 
$f(a_3)=f(a_1)*f(a_2)$. It is equivalent to $f(a_1)=f(a_3)\,\bar{*}\,f(a_2)$.
\end{definition}

The number of quandle colorings of a link diagram is a link invariant, and, more generally, considering colorings is the first step in defining many more powerful invariants, including the ones derived from quandle homology and cohomology theories; see, for example, \cite{CJKLS}.

Now we define the structure suitable for the problem of deciding whether two link diagrams with a relation on the set of components (indicating when a component can move over another component) are connected by a sequence of moves permitted by the relation.

\begin{definition} \label{parquandrel}
A {\it partial quandle with a binary relation} is a set $X$ with two partial operations 
$(x,y) \mapsto x*y$, $(x,y) \mapsto x\,\bar{*}\,y$, and a binary relation $\R\subseteq X\times X$ such that:
\begin{enumerate}
\item[(PQ1)] if $x\R y$, then $x*y$, and $x\,\bar{*}\,y$ are defined,
\item[(PQ2)] if $x\R y$, then $x*y\sim_{\R}x$, and $x\,\bar{*}\,y\sim_{\R}x$,
\item[(PQ3)] if $x\R x$, then $x*x=x=x\,\bar{*}\,x$,
\item[(PQ4)] if $x\R y$, then $x*y\,\bar{*}\,y=x=x\,\bar{*}\,y*y$,
\item[(PQ5)] if $x\R y$, $x\R z$, and $y\R z$, then $(x*y)*z=(x*z)*(y*z)$.
\end{enumerate}
\end{definition} 

\begin{definition}
A {\it partial rack with a binary relation} is a structure satisfying the axioms (PQ1), (PQ2), (PQ4), and (PQ5) of Def. \ref{parquandrel}.
\end{definition}

Just like in the case of standard racks, there are versions of the right-hand distributivity involving the second operation
$\bar{*}$, that follow from the distributivity that uses only $*$.

\begin{lemma}
Let $X$ be a structure satisfying axioms (PQ1), (PQ2), and (PQ4). Then the following conditions are equivalent:
\begin{enumerate}
\item[(1)] if $x\R y$, $x\R z$, and $y\R z$, then $(x*y)*z=(x*z)*(y*z)$,
\item[(2)] if $x\R y$, $x\R z$, and $y\R z$, then $(x\,\bar{*}\,y)*z=(x*z)\,\bar{*}\,(y*z)$,
\item[(3)] if $x\R y$, $x\R z$, and $y\R z$, then $(x*y)\,\bar{*}\,z=(x\,\bar{*}\,z)*(y\,\bar{*}\,z)$,
\item[(4)] if $x\R y$, $x\R z$, and $y\R z$, then $(x\,\bar{*}\,y)\,\bar{*}\,z=(x\,\bar{*}\,z)\,\bar{*}\,(y\,\bar{*}\,z)$.
\end{enumerate}
\end{lemma}
\begin{proof}
(1)$\Rightarrow$(2): $x\R y$, so $x\,\bar{*}\,y$ exists, and $x\,\bar{*}\,y\sim_{\R}x$; thus, 
$x\,\bar{*}\,y\R y$, $x\,\bar{*}\,y\R z$. Also $y\R z$, therefore:
\[x*z=((x\,\bar{*}\,y)*y)*z=((x\,\bar{*}\,y)*z)*(y*z).\] 
Now, $(x\,\bar{*}\,y)*z\sim_{\R}x$, $y*z\sim_{\R}y$, and $x\R y$, so $(x\,\bar{*}\,y)*z\,\R\, y*z$, thus:
\[(x*z)\,\bar{*}\,(y*z)=(((x\,\bar{*}\,y)*z)*(y*z))\,\bar{*}\,(y*z)=(x\,\bar{*}\,y)*z.\]
(2)$\Rightarrow$(1):
\[x*z=((x*y)\,\bar{*}\,y)*z=((x*y)*z)\,\bar{*}\,(y*z),\ \textrm{so}\ \]
\[(x*z)*(y*z)=(((x*y)*z)\,\bar{*}\,(y*z))*(y*z)=(x*y)*z.\]
(3)$\Leftrightarrow$(4): Proved in the same way as (1)$\Leftrightarrow$(2), replacing $*$ with $\bar{*}$, and vice versa.\\ 
(2)$\Rightarrow$(4): $x\,\bar{*}\,z\sim_{\R}x$, $y\,\bar{*}\,z\sim_{\R}y$, so
$x\,\bar{*}\,z\,\R\,y\,\bar{*}\,z$, $x\,\bar{*}\,z\,\R\,z$, and $y\,\bar{*}\,z\,\R\,z$. Thus:
\[((x\,\bar{*}\,z)\,\bar{*}\,(y\,\bar{*}\,z))*z=((x\,\bar{*}\,z)*z)\,\bar{*}\,((y\,\bar{*}\,z)*z)=x\,\bar{*}\,y,\ \textrm{and}\ \]
\[(x\,\bar{*}\,z)\,\bar{*}\,(y\,\bar{*}\,z)=(((x\,\bar{*}\,z)\,\bar{*}\,(y\,\bar{*}\,z))*z)\,\bar{*}\,z=(x\,\bar{*}\,y)\,\bar{*}\,z.\]
(4)$\Rightarrow$(2): $x*z\,\R\,y*z$, $x*z\,\R\, z$, and $y*z\,\R\, z$, thus:
\[((x*z)\,\bar{*}\,(y*z))\,\bar{*}\,z=((x*z)\,\bar{*}\,z)\,\bar{*}\,((y*z)\,\bar{*}\,z)=x\,\bar{*}\,y,\ \textrm{therefore,}\ \]
\[(x*z)\,\bar{*}\,(y*z)=(x\,\bar{*}\,y)*z.\]
\end{proof}

Now, we define two types of colorings of diagrams of links with binary relations on the set of components, using partial quandles with binary relations. For each type, the number of colorings of a diagram is an invariant; we note however, that they differ in their use: the first type seems to be more effective in distinguishing between links with relations; the second type is a basis for homological invariants. Both types utilize monotone maps between the set of arcs of a diagram, with relation induced from a relation on the set of components, to a given partial quandle with a binary relation.

\begin{definition}
Let $D$ be a diagram of a link with a binary relation on the set of its components. Let $cr$ be a crossing of $D$ such that its under-arcs belong to the component $C_i$, and its over-arc belongs to the component $C_j$. We call $cr$ a
{\it good crossing} if $C_i R C_j$, otherwise, we call it a {\it bad crossing}. A bad crossing is rigid (up to planar isotopy), in the sense that the component $C_j$ is not allowed to move over the component $C_i$. On diagrams, bad crossings will be denoted by circles around them.
\end{definition} 

\begin{figure}
\begin{center}
\includegraphics[height=11 cm]{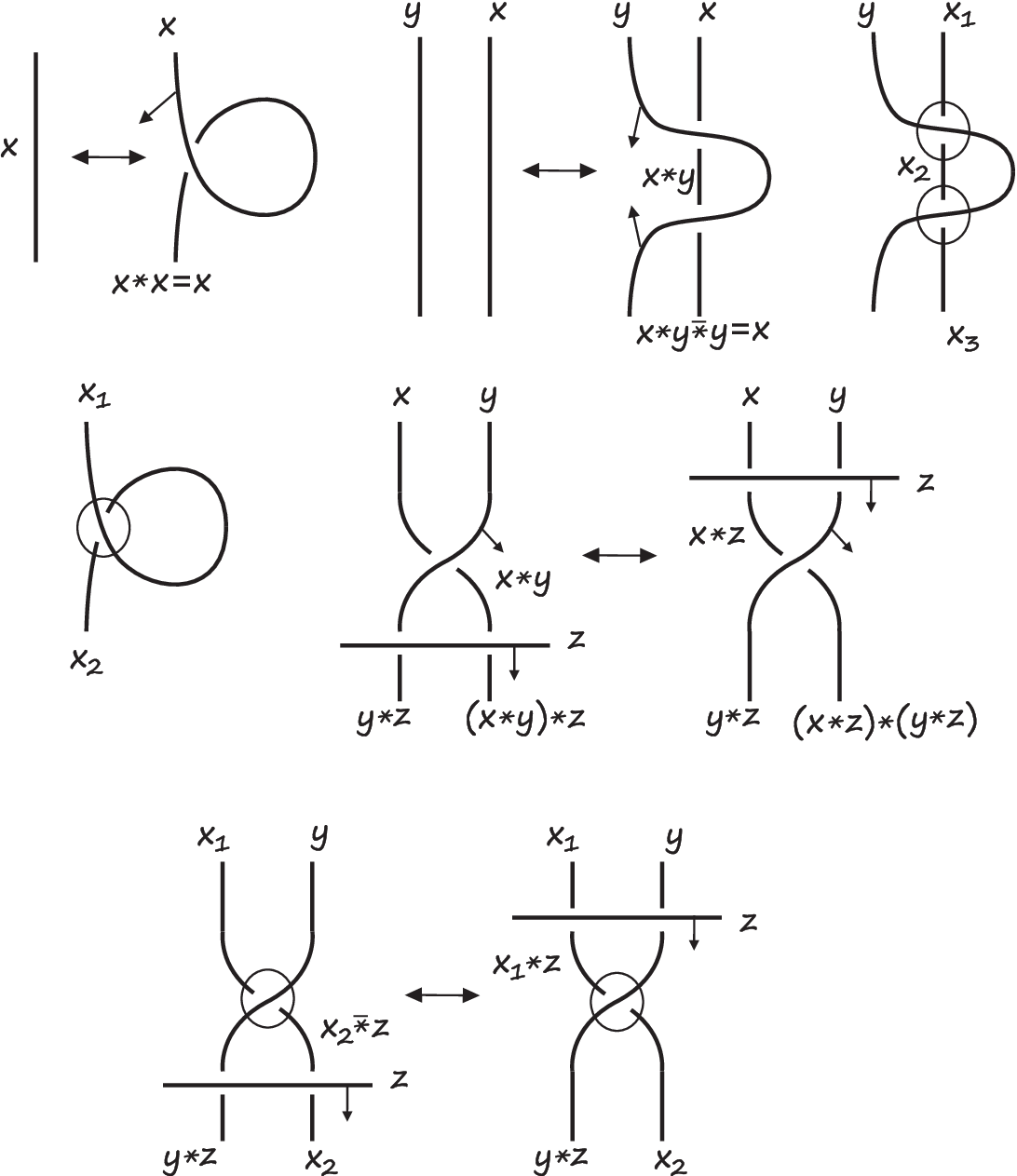}
\caption{}\label{pqcolstype1}
\end{center}
\end{figure}

\begin{definition} [Coloring of type I] \label{col1}
Let $Arcs(D)$ denote the set of arcs of an oriented link diagram with a binary relation $R$ on the set of its components,
and let $X$ be a partial quandle with a binary relation $\R$.
In this setting, a coloring of type I is a monotone map $f\colon Arcs(D)\to X$ such that the following conditions are satisfied:
\begin{enumerate}
\item[(1)] At a good crossing, the usual quandle coloring rule, depicted in Fig. \ref{qcols} is used.
\item[(2)] At a bad crossing, the quandle coloring rule is not applied; the elements assigned to the under-arcs of a bad crossing are arbitrary, as long as the requirements coming from other crossings are satisfied, and the map $f$ remains monotone.
\end{enumerate}
These conditions (for some choices of orientation) are depicted in Figure \ref{pqcolstype1}.
\end{definition}

\begin{figure}
\begin{center}
\includegraphics[height=7 cm]{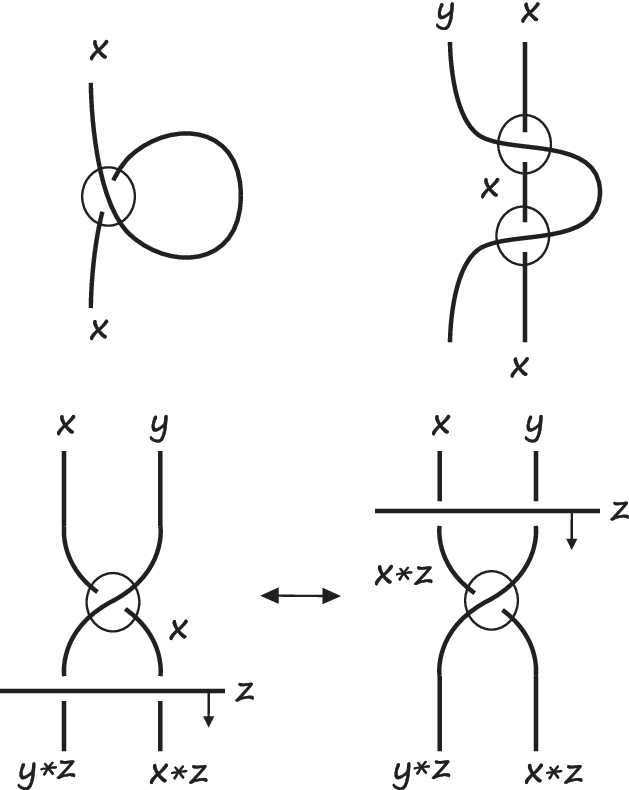}
\caption{}\label{pqcolstype2}
\end{center}
\end{figure}

\begin{definition} [Coloring of type II]
Given $Arcs(D)$, $R$, $X$, and $\R$ as in Definition \ref{col1}, we define a coloring of type II as a monotone map $f\colon Arcs(D)\to X$ satisfying:
\begin{enumerate}
\item[(1)] At a good crossing, the usual quandle coloring rule, depicted in Fig. \ref{qcols} is used.
\item[(2)] A bad crossing does not influence the coloring, that is, the colors of both under-arcs of a bad crossing are always the same.
\end{enumerate}
The condition (2) is illustrated in Fig. \ref{pqcolstype2}.
\end{definition} 

\begin{lemma}
Let $D$ be a diagram of a link with a binary relation $R$ on the set of its components, and let $X$ be a partial quandle with a binary relation. Then the number of colorings of $D$ of type I, and the number of colorings of $D$ of type II (both using $X$), are invariants under the Reidemeister moves permitted by the relation $R$.
\end{lemma}
\begin{proof}
Because the coloring maps are monotone, whenever a move allowed by the relation $R$ is performed, the partial quandle operations are defined, so that the colors of the arcs after the move can be calculated. The only exception to this is when a bad crossing occurs in the third Reidemeister move,
but just like in the case of all the other moves permitted by $R$, there is a bijection between the set of colorings before and after the move.
\end{proof}

\begin{figure}
\begin{center}
\includegraphics[height=5 cm]{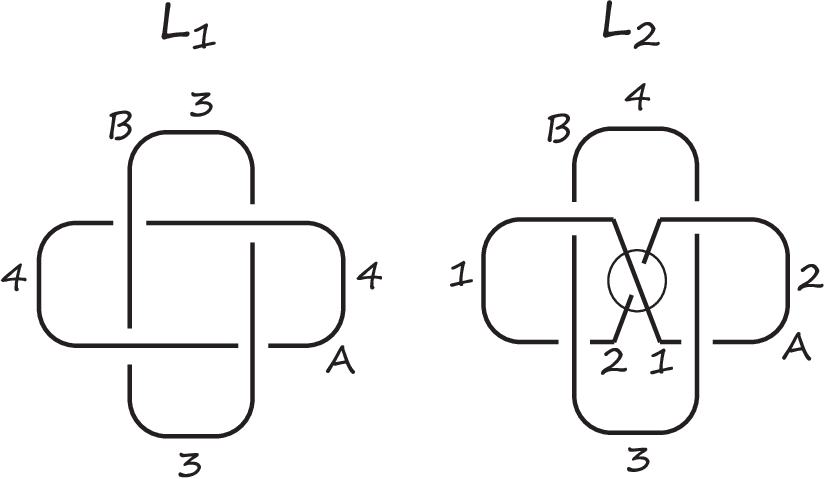}
\caption{}\label{examplcol}
\end{center}
\end{figure}

\begin{example}
Consider a partial quandle $X$ with an involutory operation $*=\bar{*}$ given by the table below, with a binary
relation: 
\[\R=\{(1,3), (1,4), (2,3), (2,4), (3,1), (3,2), (3,3), (3,4), (4,1), (4,2), (4,3), (4,4)\}.\]
\begin{center}
\begin{tabular}{c| c cccc} 
$*$& 1 & 2 & 3 & 4 \\
\hline 
1 &   &   & 2 & 1 \\
2 &   &   & 1 & 2 \\
3 & 4 & 4 & 3 & 3 \\
4 & 3 & 3 & 4 & 4 \\
\end{tabular}
\end{center}

It is not possible to obtain a distributive structure if the remaining entries in the table are filled with elements of $X$.

Let $L_1$ and $L_2$ be the two links depicted in Fig. \ref{examplcol}, with components denoted by $A$ and $B$, with the following relation $R$:

\begin{center}
\begin{tabular}{c| c cc} 
$R$& $A$ & $B$ \\
\hline 
$A$ & 0 & 1 \\
$B$ & 1 & 1 \\
\end{tabular}
\end{center}

The sets of arcs of $L_1$ and $L_2$ are equipped with the binary relations induced from $R$.
Consider the colorings (both types) of these links by $X$. Because $B R B$, and the colorings are monotone, the only colors that can be used for the components $B$ are 3 and 4. For the components $A$, one has to consider all elements of $X$.
It is not possible to color $L_1$, if the colors 1 or 2 are used for the component $A$, but there are 4 colorings of $L_1$ if both components get the colors from the set $\{3,4\}$. On the other hand, there are 8 colorings of $L_2$; one of the colorings that uses all the colors is shown in Fig. \ref{examplcol}. For these links, and this particular $X$, the colorings of type I coincide with the colorings of type II.
\end{example}

Our next goal is to define homology of partial quandles with binary relations, and obtain homological invariant of links with relations on the set of components. First we define a family of homologies of binary relations, and then we merge it with the standard rack and quandle homology by choosing appropriate chain groups, and using the usual differential. 

\subsection{Homology of binary relations}

Given a set $X$, we define a family of homologies for a binary relation $\R\subseteq X\times X$. A different homology of relations was introduced in \cite{Dow}, also, compare with the homology of reflexive and symmetric relations defined in \cite{Sos}.

\begin{definition} \label{defect}
Given an $n$-tuple $w=(x_1,\ldots,x_n)$ of elements of $X$, let $x_i$ and $x_j$ be  members of $w$ such that $i<j$. We say that there is an arrow between $x_i$ and $x_j$ if $x_i\R x_j$. We define the {\it defect} of $w$, $df(w)$, as the difference between the maximal possible number of arrows in an $n$-tuple (equal to $n(n-1)/2$) and the number of arrows of $w$. Thus, defect can be seen as the number of `missing' left-to-right arrows.
\end{definition}

\begin{figure}
\begin{center}
\includegraphics[height=1.3 cm]{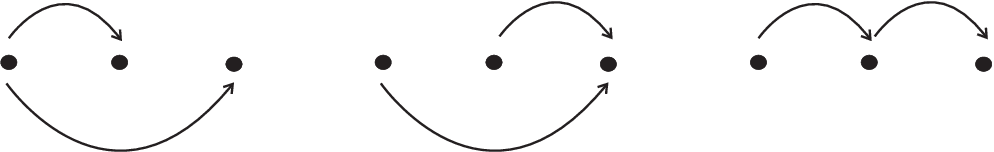}
\caption{}\label{defect1}
\end{center}
\end{figure}

Diagrams in Figure \ref{defect1} represent triples with defect 1.

\begin{definition}
For a given set $X$ with a binary relation $\R\subseteq X\times X$, let $C^k_n(X)$, for $n\geq 1$, be the free abelian 
group generated by $n$-tuples $(x_1,x_2,\ldots,x_n)$ of elements of $X$ such that $df(w)\leq k$, and let 
$C^k_0(X)=\mathbb{Z}$. 
Define a boundary   
homomorphism $\partial_n\colon C^k_n(X)\to C^k_{n-1}(X)$ by:
\[\partial_n(x_1,x_2,\ldots,x_n) = \sum_{i=1}^n (-1)^i(x_1,\ldots,x_{i-1},x_{i+1},\ldots,x_n)\] 
for $n\geq 2$, and $\partial_1=0$.
We call $(C^k_*(X),\partial_*)$ the {\it defect $k$ chain complex of the relation $\R$}, and its homology, $H^k_*(\R)$, the {\it defect $k$ homology of $\R$}.
\end{definition}

The above definition is possible because removing elements from an $n$-tuple can only decrease defect or leave it unchanged.

\subsection{Homology for partial racks with binary relations}

Rack homology and homotopy theory were first defined and 
studied in \cite{FRS}, and a modification to quandle homology  
was given in \cite{CJKLS} to define knot invariants in a state-sum form (so-called cocycle knot invariants). Later, homology of distributive structures was studied by J.H. Przytycki and his co-authors (see, for example, \cite{PS}). We recall some of these definitions, as we are going to adapt them to partial structures with binary relations. 

\begin{definition} \label{rackhom}
\begin{enumerate}
\item[(i)]
For a given rack $X$, let $C^R_n(X)$ be the free abelian 
group generated by $n$-tuples $(x_1,x_2,\ldots,x_n)$ of elements of $X$; 
in other words, $C^R_n(X) = {\Z}X^n = ({\Z}X)^{\otimes n}$, for $n\geq 1$. Let $C^R_0(X)=\mathbb{Z}$.
Define a boundary   
homomorphism $\partial_n\colon C^R_n(X) \to C^R_{n-1}(X)$ by:
\[\partial_n(x_1,x_2,\ldots,x_n) =  
\sum_{i=1}^n (-1)^i((x_1,\ldots,x_{i-1},x_{i+1},\ldots, 
x_n) -\] \[(x_1*x_i,x_2*x_i,\ldots,x_{i-1}*x_i,x_{i+1},\ldots,x_n))\] 
for $n\geq 2$, and $\partial_1=0$.
$(C^R_*(X),\partial_*)$ is called the rack chain complex of $X$.
\item[(ii)] Let $C^D_n(X)$ be a subgroup of $C^R_n(X)$ generated by $n$-tuples $(x_1,\ldots,x_n)$ with $x_i=x_{i+1}$ for some $i\in\{1,\ldots,n-1\}$, if $n\geq 2$, and let $C^D_n(X)=0$ otherwise. If $X$ is a quandle, then
$(C^D_*(X),\partial_*)$ is a subchain 
complex of $(C^R_*(X),\partial_*)$, called the degenerate chain complex of a quandle $X$.
\item[(iii)] The quotient chain complex $(C^Q_*(X),\partial_*)$, obtained by taking $C^Q_n(X)=C^R_n(X)/C^D_n(X)$, is 
called the quandle chain complex. 
\item[(iv)] The homology  of the rack, degenerate, and quandle chain complexes 
is called the rack, degenerate, and quandle homology, and is denoted by $H^R_*(X)$, $H^D_*(X)$, and $H^Q_*(X)$, respectively.
\end{enumerate}
\end{definition}

Now we incorporate the notion of defect and partial operation $*$ into the above homologies.

\begin{theorem}
Let $X$ be a partial rack with a binary relation $\R$. Let $C^{PR}_n(X)$, for $n\geq 1$, be the free abelian 
group generated by $n$-tuples $w=(x_1,x_2,\ldots,x_n)$ of elements of $X$ such that $df(w)=0$; also,
let $C^{PR}_0(X)=\mathbb{Z}$. The boundary $\partial_n\colon C^{PR}_n(X)\to C^{PR}_{n-1}(X)$ is defined as in Def. \ref{rackhom}. Then $(C^{PR}_*(X),\partial_*)$ is a chain complex. Let $C^{PD}_n(X)$ be a subgroup of $C^{PR}_n(X)$ generated by $n$-tuples $(x_1,\ldots,x_n)$ with defect 0, with $x_i=x_{i+1}$ for some $i\in\{1,\ldots,n-1\}$, if $n\geq 2$, and let $C^{PD}_n(X)=0$ otherwise. If $X$ is a partial quandle with a binary relation, then $(C^{PD}_*(X),\partial_*)$ is a subchain 
complex of $(C^{PR}_*(X),\partial_*)$, and we can form a quotient chain complex $(C^{PQ}_*(X),\partial_*)$
by taking $C^{PQ}_n(X)=C^{PR}_n(X)/C^{PD}_n(X)$. The homology of this last complex, $H^{PQ}_*(X)$, will be called the homology of a partial quandle with a binary relation, and the homology of $(C^{PR}_*(X),\partial_*)$, denoted 
$H^{PR}_*(X)$, will be called the homology of a partial rack with a binary relation.
\end{theorem}
\begin{proof}
The fact that we use defect 0 in the above definition means that in the $n$-tuple $(x_1,\ldots,x_n)\in C^{PR}_n(X)$,
$x_i\R x_j$, for every $i<j$, and thus, the multiplication $x_i*x_j$ is defined. Moreover, from the definition of 
a partial rack with a binary relation follows that $x_i*x_j\sim_{\R}x_i$. Thus, after removing an element $x_j$ from the $n$-tuple, and multiplying all the preceding elements by it, the defect of the new $(n-1)$-tuple remains 0. As we noted when defining the homology of a binary relation, simply removing an element from an $n$-tuple (as is the case in the first part of the differential) cannot increase the defect. Thus, $\partial_n$ indeed takes $C^{PR}_n(X)$ to $C^{PR}_{n-1}(X)$. 
Because the formula for $\partial_n$ is the same as in standard rack homology, and all the required multiplications are defined, it follows that $\partial_{n-1}\circ\partial_{n}=0$.

Let $X$ be a partial quandle with a binary relation. Then, if an $n$-tuple $(x_1,\ldots,x_i,x_i,\ldots,x_n)$ has defect 0, it follows, in particular, that $x_i\R x_i$, and $x_i*x_i=x_i$; also $x_i*x_j\sim_{\R}x_i$, for any $j>i$, so $x_i*x_j\,\R\, x_i*x_j$.
Then it follows, as in the standard quandle homology, that $(C^{PD}_*(X),\partial_*)$ is a subchain 
complex of $(C^{PR}_*(X),\partial_*)$, that is: $\partial_n(C^{PD}_n(X))\subseteq C^{PD}_{n-1}(X)$.
\end{proof} 

\begin{remark}
If $X$ is a set with a relation $\R\subseteq X\times X$, and a binary operation $(x,y) \mapsto x*y$, defined for all $x$, $y\in X$, and satisfying the conditions:
\begin{enumerate}
\item[(1)] $x\leadsto_{\ \ {\R}} x*y$, for all $x$, $y\in X$,
\item[(2)] $(x*y)*z=(x*z)*(y*z)$, for all $x$, $y$, $z\in X$,
\end{enumerate}
then a homology of such a structure can be defined for any defect. Indeed, let $C^k_n(X)$, for $n\geq 1$, be the free abelian group generated by $n$-tuples $(x_1,x_2,\ldots,x_n)$ of elements of $X$ such that $df(w)\leq k$, and let $C^k_0(X)=\mathbb{Z}$. 
Define a boundary homomorphism $\partial_n: C^k_n(X) \to C^k_{n-1}(X)$ as in Definition \ref{rackhom}.
Then the condition (1) assures that $\partial_n(C^k_n(X))\subseteq C^k_{n-1}(X)$, and the condition (2)
gives $\partial_{n-1}\circ\partial_{n}=0$. The notion of defect gives a filtration, and the corresponding spectral sequence should be investigated. 
\end{remark}   

Let us recall the procedure of assigning a 2-cycle in rack homology
to an oriented colored link diagram, introduced in \cite{Gre}. Let a link diagram be colored with elements of a rack $X$ according to the rules depicted in Fig. \ref{qcols}.
Each positive crossing represents a pair $(x,y)\in C^R_2(X)$, where $x=f(a_1)$ is the color of the under-arc away 
from which points the normal of the over-arc labeled $y=f(a_2)$. In the case of a negative crossing, we write $-(x,y)$.
The sum of such 2-chains taken over all crossings of the diagram forms a 2-cycle. Thus, it represents an element in $H^R_2(X)$.

We can adjust this procedure to work with homology of partial quandles with binary relations. 

\begin{theorem}
Let $X$ be a partial quandle with a binary relation $\R$, and consider a coloring $f\colon Arcs(D)\to X$ of type II
of an oriented link diagram $D$ with a binary relation $R$ on the set of its components $C=\{C_1,\ldots,C_s\}$.
Let $Gcr_i$, $i\in\{1,\ldots,s\}$, be the set of all good crossings in which the under-arcs belong to the component $C_i$.
To each positive good crossing assign the pair $(x,y)\in X\times X$, where $x$ is the color of the under-arc away 
from which points the normal of the over-arc with color $y$; to a negative crossing, assign $-(x,y)$.
Now, let $z_i$, $i\in\{1,\ldots,s\}$, be the sum of such signed pairs of colors, taken over all crossings from $Gcr_i$.
Then each $z_i$ is a cycle in $(C^{PQ}_2(X),\partial_2)$, and it represents an element of $H^{PQ}_2(X)$ that is invariant under all the Reidemeister moves permitted by the relation $R$.
\end{theorem}
\begin{proof}
$z_i$ being a cycle follows from the fact that at a good crossing the pair $(x,y)$ has defect 0, $x*y$ is defined, and 
\[\partial_2(x,y)=-y+y+x-x*y=x-x*y,\] that is, only the colors of the under-arcs appear in the boundary of a pair of colors assigned to a crossing. When traveling along the component $C_i$, ignoring bad crossings, and writing the boundaries of signed pairs of colors assigned to crossings from $Gcr_i$, reductions occur for each pair of consecutive crossings, and since the component is closed, the sum of boundaries is zero.

When a kink is created or deleted as a result of a first Reidemeister move on the component $C_i$ such that $C_i R C_i$, then the pair of colors $(x,x)$ involved in the move has defect 0 (due to monotonicity of colorings), and so it belongs to 
$C^{PD}_2(X)$, and represents 0 in $H^{PQ}_2(X)$. 

A second Reidemeister move in which two crossings belonging to $Gcr_i$ are created or deleted, only adds or deletes
$(x,y)-(x,y)=0$.

In a third Reidemeister move involving a bad crossing, the two pairs of colors corresponding to the good crossings are not changed by the move, because of the way bad crossings are colored in a coloring of type II.

Finally, the third Reidemeister move in which all the crossings are good and positive, and in which the lowest component is $C_i$, adds 
\[\pm \partial_3(x,y,z)=-(y,z)+(y,z)+(x,z)-(x*y,z)-(x,y)+(x*z,y*z)\] to $z_i$, thus, the homology class of $z_i$ is not changed (see Fig. \ref{pqcolstype1} for this move). The third Reidemeister move in which all the crossings are positive, together with the second Reidemeister moves, generates all the other third Reidemeister moves with different orientations.
Moreover, by looking at the proofs in \cite{Pol} describing generating sets of moves, we notice that all the versions of the third Reidemeister move can be generated so that the bottom component never moves over the middle or the top component, and the middle component never moves over the top component; it means that all the third Reidemeister moves can be generated in our conditional setting, so the proof is ended. 
\end{proof}

\begin{figure}
\begin{center}
\includegraphics[height=5 cm]{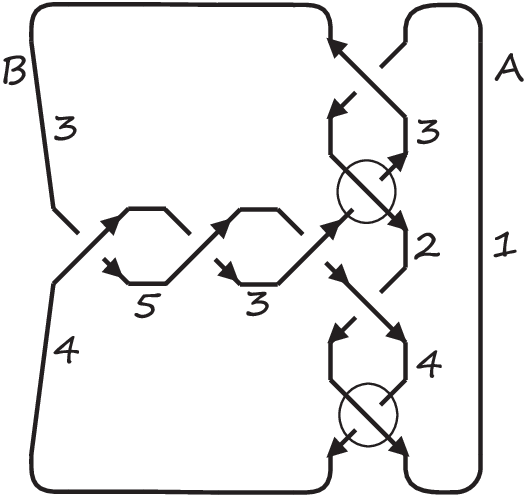}
\caption{}\label{gentor3}
\end{center}
\end{figure}

\begin{example}
Consider the following involutory quandle $X$, which satisfies the axioms of a partial quandle with a binary relation 
\[ \R=\{(1,1),(1,2),(1,3),(1,4),(1,5),(2,1),(2,2),(2,3),(2,4),\]
\[(2,5),(3,3),(3,4),(3,5),(4,3),(4,4),(4,5),(5,3),(5,4),(5,5)\}.\]

\begin{center}
\begin{tabular}{c| c ccccc} 
$*$& 1 & 2 & 3 & 4 & 5 \\
\hline 
1 & 1 & 1 & 2 & 2 & 2 \\
2 & 2 & 2 & 1 & 1 & 1 \\
3 & 3 & 3 & 3 & 5 & 4 \\
4 & 4 & 4 & 5 & 4 & 3 \\
5 & 5 & 5 & 4 & 3 & 5 \\
\end{tabular}
\end{center}

The calculations in \cite{GAP4} show that the torsion parts of the homology groups $H_i^{PQ}(X)$, for $i=2,\ldots,7$, are as follows: $\Z_3$, $\Z_3^3$, $\Z_3^5$, $\Z_3^8$, $\Z_3^{13}$, $\Z_3^{20}$. 
Compare it with the torsion parts in standard quandle homology $H_i^{Q}(X)$, $i=2,\ldots,5$:
$\Z_3$, $\Z_3^4$, $\Z_3^{10}$, $\Z_3^{23}$. The calculations for $H_i^{PQ}(X)$ are much faster due to the smaller size of the chain groups. 

The link depicted in Fig. \ref{gentor3} is assumed to have the following relation on the set of components: 
\begin{center}
\begin{tabular}{c| c cc} 
$R$& $A$ & $B$ \\
\hline 
$A$ & 1 & 1 \\
$B$ & 0 & 1 \\
\end{tabular}
\end{center}

Then, Fig. \ref{gentor3} shows a coloring of type II, using $X$, that represents a torsion $\Z_3$ in $H_2^{PQ}(X)$.
The cycle $c\in C_2^{PQ}(X)$ corresponding to this coloring is: 
\[-(5,4)-(3,5)-(4,3)-(1,4)+(1,3).\]
It is a sum of two cycles yielded by the components $A$ and $B$: $c=c_1+c_2$, where
$c_1=-(1,4)+(1,3)$, and $c_2=-(5,4)-(3,5)-(4,3)$. $c_1$ is the part that represents torsion, and $c_2$ is a boundary:
\[c_2=\partial_3(-(3,4,5)-(5,3,4)+(1,2,1)+(3,4,3)).\]
\end{example}

\end{document}